\documentclass[11pt]{article}

\usepackage{color}
\usepackage{latexsym}
\usepackage{amssymb}
\usepackage{graphicx}
\usepackage{amsmath, amsfonts,amssymb,theorem,euscript,
array,enumerate,amsfonts,mathrsfs}

\newtheorem{Theorem}{Theorem}[part]

\newtheorem{Proposition}{Proposition}[part]

\newtheorem{Lemma}{Lemma}[part]

\newtheorem{Remark}{Remark}[part]

\def \Frac{\displaystyle\frac}

\def \b1{\bf{1}}

\def \I{\mathbb{I}}
\def \N{\mathbb{N}}
\def \R{\mathbb{R}}

\def \Z{\mathbb{Z}}
\def \E{\mathbb{E}}
\def \F{\mathbb{F}}
\def \P{\mathbb{P}}

\def \T{\mathbb{T}}
\def \X{\mathbb{X}}

\def\esssup_#1{\underset{#1}{\mathrm{ess\,sup\, }}}

\def \Ac{{\cal A}}

\def \Fc{{\cal F}}

\def \Mc{{\cal M}}
\def \Nc{{\cal N}}

\def \Sc{{\cal S}}

\def \ni{\noindent}

\def \eps{\varepsilon}

\def \ep{\hbox{ }\hfill$\Box$}

\def\Dt#1{\Frac{\partial #1}{\partial t}}

\def\reff#1{{\rm(\ref{#1})}}

\def\beqs{\begin{eqnarray*}}
\def\enqs{\end{eqnarray*}}
\def\beq{\begin{eqnarray}}
\def\enq{\end{eqnarray}}

\begin{document}

\title{Time discretization and quantization methods for \\  optimal multiple switching problem\thanks{We would like to thank Damien Lamberton, Nicolas Langren\'e and Gilles Pag\`es for helpful remarks.}}
\author{Paul Gassiat$^{1)}$, Idris Kharroubi$^{2)}$, Huy\^en Pham$^{1),3)}$}

\date{September 23, 2011 \\ Revised version: February  7, 2012}

\maketitle {\noindent
\small
\begin{tabular}{llll}
$^{1)}$ & Laboratoire de Probabilit\'es et              &  $^{2)}$ & CEREMADE, CNRS, UMR 7534   \\
        & Mod\`eles Al\'eatoires,  CNRS, UMR 7599              &          &  Universit\'e Paris Dauphine\\
        &  Universit\'e Paris 7 Diderot,     &          &   kharroubi at ceremade.dauphine.fr \\
        &      gassiat, pham  at math.univ-paris-diderot.fr                    &  &                   \\
$^{3)}$ &          CREST-ENSAE                    &           &   \\
        &             and Institut Universitaire de France                                  &          &    \\[1ex]
\end{tabular} }

\begin{abstract}
In this paper, we study probabilistic numerical methods based on optimal quantization algorithms for computing the solution to 
optimal multiple switching problems with regime-dependent state process. We first consider a discrete-time approximation of the  
optimal switching problem, and analyze its rate of convergence. 
Given a time step $h$, the error is in general of order $(h \log(1/h))^{1/2}$, and of order 
$h^{1/2}$ when the switching costs do not depend on the state process. We next propose quantization numerical schemes for the space 
discretization of the discrete-time Euler state process.  A  Markovian quantization approach relying on the optimal quantization of the normal distribution arising in the Euler scheme is analyzed.  In the particular case of uncontrolled state process, we  describe an alternative marginal quantization method, which extends the recursive algorithm for optimal stopping problems as in \cite{balpag03b}.  
A priori $L^p$-error estimates 
are stated in terms of quantization errors. Finally, some numerical tests are performed for an optimal switching problem with two regimes. 
\end{abstract}

\vspace{5mm}


\vspace{5mm}


\newpage

\section{Introduction}

\setcounter{equation}{0} \setcounter{Assumption}{0}
\setcounter{Theorem}{0} \setcounter{Proposition}{0}
\setcounter{Corollary}{0} \setcounter{Lemma}{0}
\setcounter{Definition}{0} \setcounter{Remark}{0}

On some filtered probability space $(\Omega,\Fc,\F=(\Fc_t)_{_{t\geq 0}},\P)$, let us introduce the controlled regime-switching diffusion in $\R^d$ governed by
\beqs
dX_t &=& b(X_t,\alpha_t) dt + \sigma(X_t,\alpha_t) dW_t,
\enqs
where $W$ is a standard $d$-dimensional Brownian motion, 
$\alpha$ $=$ $(\tau_n,\iota_n)_n$ $\in$ $\Ac$ is the switching control represented by a nondecreasing sequence of stopping times 
$(\tau_n)$ together with a sequence $(\iota_n)$ of $\Fc_{\tau_n}$-measurable random variables valued in a finite set $\{1,\ldots,q\}$, and 
$\alpha_t$ is the current regime process, i.e. 
$\alpha_t$ $=$ $\iota_n$ for $\tau_n\leq t <\tau_{n+1}$.  We then consider the optimal switching problem 
over a finite horizon:
\beq \label{defJ}
V_0 &=& \sup_{\alpha\in\Ac}   \E \Big[ \int_0^T f(X_t,\alpha_t) dt + g(X_T,\alpha_T) - \sum_{\tau_n\leq T} c(X_{\tau_n},\iota_{n-1},\iota_n) \Big].  
\enq
Optimal switching problems can be seen as sequential optimal stopping problems belonging to the class of impulse control problems, and arise in many applied fields, for example in real option pricing in economics and finance. It has attracted a lot of interest during the past decades, and we refer to Chapter 5 in the book \cite{pha09}  and the references therein for a survey of  some applications and results in this topic. It is well-known that optimal switching problems are related via the dynamic programming approach to a system of variational inequalities with inter-connected obstacles in the form:
\beq
 \min \Big[ - \Dt{v_i} - b(x,i).D_x v_i - \frac{1}{2} {\rm tr}(\sigma(x,i)\sigma(x,i)'D_x^2 v_i) - f(x,i) \; , \;  & & \label{sysvar} \\
v_i - \max_{j\neq i}(v_j - c(x,i,j)) \Big]  &=& 0 \;\; \mbox{ on } [0,T)\times\R^d,  \nonumber
\enq
together with the terminal condition $v_i(T,x)$ $=$ $g(x,i)$, for any $i$ $=$ $1,\ldots,q$.  Here $v_i(t,x)$ is the value function to the 
optimal switching problem starting at time $t$ $\in$ $[0,T]$  from the state $X_t$ $=$ $x$ $\in$ $\R^d$ and  the regime $\alpha_t$ $=$ $i$ $\in$ $\{1,\ldots,q\}$, and the solution to the system \reff{sysvar} has to be understood in the weak sense, e.g. viscosity sense.  

The purpose of this paper is to solve numerically the optimal switching problem \reff{defJ}, and consequently the system of variational inequalities  \reff{sysvar}.   These equations can be solved by analytical methods (finite differences, finite elements, etc ...), see e.g. \cite{mar05}, but are known to require heavy computations, especially in high dimension. Alternatively, when the state process is uncontrolled, i.e. regime-independent,  optimal switching problems are connected to multi-dimensional reflected Backward Stochastic Differential Equations (BSDEs) with oblique reflections, as shown in \cite{hamzha10} and \cite{hutan10}, and the recent paper \cite{chaelikha10} introduced a discretely obliquely reflected numerical scheme to solve such BSDEs.  From a computational viewpoint, there are rather few papers dealing with numerical experiments for optimal switching problems.  The special case of two regimes for switching problems  can be reduced to the re\-solution of a single BSDE with two reflecting barriers when considering the difference value process, and is exploited numerically in \cite{hamjea07}. 
We mention also the paper \cite{carlud08}, which solves an optimal switching problem with three regimes by considering a cascade of reflected BSDEs with one reflecting barrier 
derived  from an iteration on the number of switches.

We propose probabilistic numerical methods based on dynamic programming  and  optimal quantization methods combined with a suitable time discretization procedure for computing  the solution to optimal multiple  switching problem.  Quantization methods were introduced in \cite{balpag03b} for solving variational inequality 
with given obstacle associated to optimal stopping problem of some diffusion process $(X_t)$. The basic idea is the following. One first  approximates  the (continuous-time) optimal stopping problem by the Snell envelope for the Markov chain  $(\bar X_{t_k})$  defined as the Euler scheme of the (uncontrolled) diffusion $X$, and then spatially discretize each random 
vector $\bar X_{t_k}$ by a random vector  taking finite values through a quantization procedure. More precisely, $(\bar X_{t_k})_k$ is approximated  by $(\hat X_k)_k$ where $\hat X_k$ is the projection of $\bar X_{t_k}$ 
on a finite grid   in the state space following the closest neighbor rule. 
The induced $L^p$-quantization error, $\|\bar X_{t_k} - \hat X_k\|_p$, depends only on the distribution of $\bar X_{t_k}$ and the grid, which may be    chosen in order to minimize the quantization error. Such an optimal choice, called optimal quantization, is achieved  by the competitive learning vector quantization algorithm (or Kohonen algorithm) developed in full details in \cite{balpag03b}.  One finally computes the approximation of the 
optimal stopping problem by a quantization tree algorithm, which mimics the backward dynamic programming of the Snell envelope. 
In this paper, we develop quantization methods  to our general framework of optimal switching problem.  With respect to standard optimal stopping problems,  some new features arise  
on one hand from  the regime-dependent state process,  and on the other hand from the multiple switching times, and the discrete sum for the cumulated  switching costs.

We first study a time discretization of the optimal switching problem by considering an Euler-type scheme  with step $h$ $=$ $T/m$ for the regime-dependent state process $(X_t)$  controlled by the switching strategy $\alpha$:
\beq \label{eulertype}
\bar X_{t_{k+1}} &=& \bar X_{t_k} + b(\bar X_{t_k},\alpha_{t_k}) h  + \sigma(\bar X_{t_k},\alpha_{t_k})\sqrt{h} \; \vartheta_{k+1}, \;\;
t_k = kh, \; k =0,\ldots,m,
\enq
where $\vartheta_k$, $k$ $=$ $1,\ldots,m$, are iid, and $\Nc(0,I_d)$-distributed. We then introduce the optimal switching problem for the discrete-time process $(\bar X_{t_k})$ controlled by  switching strategies with stopping times valued in the discrete time grid $\{t_k, k = 0,\ldots,m\}$.  The convergence of this  discrete-time problem is analyzed, and we prove that the error is in general of order $(h \log(1/h))^{\frac{1}{2}}$, and of order $h^{\frac{1}{2}}$,
as for optimal stopping problems, when the switching costs $c(x,i,j) \equiv c(i,j)$ do not depend on the state process. 
Arguments of the proof rely on a regularity result of the controlled diffusion with respect to the switching strategy, and moment estimates on the number of switches. 
This improves and extends the convergence rate  result  in \cite{chaelikha10} derived in the case where $X$ is regime-independent. 

Next, we propose approximation schemes by quantization for computing explicitly 
the solution to the discrete-time optimal switching problem.  Since the controlled Markov chain  $(\bar X_{t_k})_k$ 
cannot be directly quantized as in standard optimal stopping problems, we adopt a Markovian quantization approach in the spirit of \cite{pagphapri04}, by  considering an  optimal 
quantization of the Gaussian random vector $\vartheta_{k+1}$  arising in the Euler scheme \reff{eulertype}. A quantization tree algorithm is then designed for computing the approximating 
value function, and we provide error estimates in terms of the quantization errors $\|\vartheta_k - \hat\vartheta_k\|_p$ and state space grid parameters.  Alternatively, in the case 
of regime-independent state process,  we propose a quantization algorithm in the vein of \cite{balpag03b} based on marginal quantization of the uncontrolled Markov chain 
$(\bar X_{t_k})_k$.  A priori $L^p$-error estimates are also established in terms of quantization errors $\| \bar X_{t_k}-\hat X_k\|_p$.  Finally, some numerical tests on the two 
quantization algorithms are performed for an optimal switching problem with two regimes.

The plan of this paper is organized as follows.  Section 2 formulates the optimal swit\-ching problem and sets the standing  assumptions.  We also show some preliminary results about moment estimates on  the number of switches. We describe in Section 3  the time discretization procedure, and study the rate of convergence of the discrete-time approximation for the optimal switching problem. Section 4 is devoted to the approximation schemes by quantization for the explicit computation of the value function to the discrete-time optimal switching problem, and to the error analysis.  Finally, we illustrate our results with some numerical tests in Section 5.

\section{Optimal  switching problem}

\setcounter{equation}{0} \setcounter{Assumption}{0}
\setcounter{Theorem}{0} \setcounter{Proposition}{0}
\setcounter{Corollary}{0} \setcounter{Lemma}{0}
\setcounter{Definition}{0} \setcounter{Remark}{0}

\subsection{Formulation and assumptions}

We formulate the finite horizon multiple switching problem.  Let us fix a finite time $T$ $\in$ $(0,\infty)$, and some filtered 
pro\-bability space $(\Omega,\Fc,\F=(\Fc_t)_{_{t\geq 0}},\P)$ satisfying the usual conditions. 
Let $\I_q$ $=$ $\{1,\ldots,q\}$ be the set of  all possible regimes (or activity modes).  A switching 
control is a double sequence $\alpha$ $=$ $(\tau_n,\iota_n)_{n\geq 0}$, where $(\tau_n)$ is a nondecreasing sequence 
of stopping times, and $\iota_n$ are $\Fc_{\tau_n}$-measurable random variables valued in $\I_q$.  The switching control 
$\alpha$ $=$ $(\tau_n,\iota_n)$ is said to be admissible, and denoted by $\alpha$ $\in$ $\Ac$, if  there exists an integer-valued random variable $N$ with $\tau_N$ $>$ $T$ a.s. 
Given $\alpha$ $=$ $(\tau_n,\iota_n)_{n\geq 0}$ $\in$ $\Ac$,  we may then associate the indicator of the regime value defined at any time $t$ $\in$ $[0,T]$ by
\beqs
I_t &=& \iota_0  {\bf 1}_{\{0 \leq t < \tau_0\}} +  \sum_{n \geq 0} \iota_n {\bf 1}_{\{\tau_n\leq t< \tau_{n+1}\}},
\enqs
which we shall sometimes identify with the switching control $\alpha$, and we introduce $N(\alpha)$  the (random) number of switches before $T$:
\beqs
N(\alpha) & = & \#\big\{n\geq 1~: \tau_n \leq T\big\}\;. 
\enqs
For $\alpha$ $\in$ $\Ac$, we consider the controlled regime-switching diffusion process valued in $\R^d$, governed by the dynamics
\beq \label{dynX}
dX_s &=& b(X_s,I_s) ds + \sigma(X_s,I_s) dW_s,  \;\;\; X_0 = x_0 \in \R^d,
\enq
where $W$ is a standard $d$-dimensional Brownian motion on $(\Omega,\Fc,\F=(\Fc_t)_{0\leq t\leq  T},\P)$. We shall assume that the coefficients  $b_i$ $=$ $b(.,i)$: $\R^d$ $\rightarrow$ $\R^d$, and $\sigma_i(.)$ $=$ $\sigma(.,i)$ $:$ $\R^d$ $\rightarrow$ $\R^{d\times d}$, $i$ $\in$ $\I_q$, satisfy the usual Lipschitz conditions. 
 
We are given a running reward, terminal gain functions $f,g$ $:$ $\R^d\times\I_q$ $\rightarrow$ $\R$, and a cost function $c$ $:$ $\R^d\times\I_q\times\I_q$ $\rightarrow$ $\R$, and we  
set $f_i(.)$ $=$ $f(.,i)$, $g_i(.)$ $=$ $g(.,i)$, $c_{ij}(.)$ $=$ $c(.,i,j)$, $i,j$ $\in$ $\I_q$.  We shall assume the Lipschitz condition: 

\vspace{2mm}

{\bf (Hl)} \hspace{2mm} The coefficients $f_i$, $g_i$ and $c_{ij}$, $i,j$ $\in$ $\I_q$ are Lipschitz continuous on $\R^d$.  

\vspace{2mm}

We also make the natural triangular condition on the functions $c_{ij}$ representing the instantaneous cost for switching from regime $i$ to $j$: 

\vspace{2mm}

{\bf (Hc)}
\beqs
c_{ii}(.) &=& 0,  \;\;\; i \in \I_q, \\
\inf_{x\in\R^d}c_{ij}(x) & > & 0, \;\;\; \mbox{ for } i,j \in \I_q,  \; j \neq i,\\
\inf_{x\in\R^d} \big[ c_{ij}(x) + c_{jk}(x) - c_{ik}(x)] &>& 0,  \;\;\; \mbox{ for } i,j,k \in \I_q,  \; j \neq i,k.
\enqs 
The triangular condition on the switching costs $c_{ij}$ in {\bf(Hc)} means that when one changes  from regime $i$ to some regime $j$, then it is not optimal to switch again immediately to another regime, since it would induce a higher total cost, and so one should stay for a while in the regime $j$.

\vspace{2mm}

The expected total profit over $[0,T]$ 
for running the system with the admissible switching control $\alpha$ $=$ $(\tau_n,\iota_n)$ $\in$ $\Ac$ is given by: 
\beqs
J_0(\alpha) &=& \E \Big[ \int_0^T f(X_t,I_t) dt + g(X_T,I_T) - \sum_{n= 1}^{N(\alpha)} c(X_{\tau_n},\iota_{n-1},\iota_n) \Big]. 
\enqs
The maximal profit is then defined  by 
\beq \label{defV0}
V_0 &=& \sup_{\alpha\in\Ac} J_0(\alpha). 
\enq
The dynamic version of this optimal switching problem is formulated as follows. 
For $(t,i)$ $\in$ $[0,T]\times\I_q$, we denote by $\Ac_{t,i}$ the set of admissible switching controls $\alpha$ $=$ $(\tau_n,\iota_n)$ starting from $i$ at time $t$, i.e. $\tau_0$ $=$ $t$, $\iota_0$ $=$ $i$.  Given $\alpha$ $\in$ $\Ac_{t,i}$, and $x$ $\in$ $\R^d$, and under the Lipschitz conditions on 
$b$, $\sigma$, there exists a unique strong solution to \reff{dynX} starting from $x$ at time $t$, and denoted by $\{X_s^{t,x,\alpha},t\leq s\leq T\}$. 
It is then given by
\beq \label{Xxalpha}
X_s^{t,x,\alpha} &=& x +  \sum_{\tau_n\leq s}  \int_{\tau_n}^{\tau_{n+1}\wedge s} b_{\iota_n}(X_u^{t,x,\alpha}) du 
+  \int_{\tau_n}^{\tau_{n+1}\wedge s}  \sigma_{\iota_n}(X_u^{t,x,\alpha}) dW_u, \; t \leq s \leq T.\qquad  
\enq
The value function of the optimal switching problem is defined by
\beq \label{defvi}
v_i(t,x) &=& \sup_{\alpha\in\Ac_{t,i}}   \E \Big[ \int_t^T f(X_s^{t,x,\alpha},I_s) ds + g(X_T^{t,x,\alpha},I_T) 
- \sum_{n= 1}^{N(\alpha)} c(X_{\tau_n}^{t,x,\alpha},\iota_{n-1},\iota_n) \Big],\qquad
\enq
for any $(t,x,i)$ $\in$ $[0,T]\times\R^d\times\I_q$, so that $V_0$ $=$ $\max_{i\in\I_q} v_i(0,x_0)$.  

For simplicity, we shall also make the assumption
\beq \label{hypg}
g_i(x) & \geq &  \max_{j\in\I_q} [ g_{j}(x) -  c_{ij}(x)], \;\;\; \forall (x,i) \in \R^d\times\I_q.
\enq
This means that any switching decision at  horizon  $T$ induces a terminal profit, which is smaller than a no-decision at this time, and is thus 
suboptimal. Therefore, the terminal condition for the value function is given by: 
\beqs
v_i(T,x) &=& g_i(x), \;\;\; (x,i) \in \R^d\times\I_q. 
\enqs
Otherwise, it is given in general by $v_i(T,x)$ $=$ $\max_{j\in\I_q} [ g_{j}(x) -  c_{ij}(x)]$.


\vspace{2mm}

\noindent {\bf Notations.}  $|.|$ will denote the canonical Euclidian norm on $\R^d$, and $(.|.)$ the corresponding inner product. 
For any $p$ $\geq$ $1$,  and $Y$ random variable on $(\Omega,\Fc,\P)$, we denote by $\| Y \|_p$ $=$ $(\E|Y|^p)^{1\over p}$.

\subsection{Preliminaries}

We first show that one can restrict the optimal switching problem to controls  $\alpha$  with bounded moments of $N(\alpha)$.   
More precisely, let us associate to a strategy $\alpha$ $\in$ $\Ac_{t,i}$,  the cumulated cost process  $C^{t,x,\alpha}$ defined by
\beqs
C^{t,x,\alpha}_u & = & \sum_{n\geq 1} c(X_{\tau_n}^{t,x,\alpha},\iota_{n-1},\iota_n) \mathbf{1}_{\tau_n\leq u}, \;\;\;\;\;\;\;   t \leq u \leq T. 
\enqs
We  then consider for $x\in\R^d$ and  $K$ $>$ $0$  the subset $\Ac^K_{t,i}(x)$ of  $\Ac_{t,i}$ defined by
 \beqs
\Ac_{t,i}^K(x) &=& \Big\{ \alpha \in \Ac_{t,i}:  \E\big|C^{t,x,\alpha}_T\big|^2  \; \leq \;  K(1+|x|^2) \Big\}.
 \enqs


\begin{Proposition}\label{proprestalpha}
Assume that {\bf (Hl)} and {\bf (Hc)} holds.  Then, there exists some positive constant $K$ s.t. 
\beq \label{relviK}
v_i(t,x) &=&  \sup_{\tiny{\alpha\in\Ac^{K}_{t,i}(x)}}   \E \Big[ \int_t^T f(X_s^{t,x,\alpha},I_s) ds + g(X_T^{t,x,\alpha},I_T)  - \sum_{n= 1}^{N(\alpha)} 
c(X_{\tau_n}^{t,x,\alpha},\iota_{n-1},\iota_n) \Big]
\enq
for any $(t,x,i)$ $\in$ $[0,T]\times\R^d\times\I_q$.
\end{Proposition}

\begin{Remark} \label{remNal}
{\rm  Under the uniformly strict positive condition on the switching costs in {\bf (Hc)}, there exists some positive constant $\eta$ $>$ $0$ s.t.  $N(\alpha)$ $\leq$ $\eta C_T^{t,x,\alpha}$ for any 
$(t,x,i)$ $\in$ $[0,T]\times\R^d\times\I_q$, $\alpha$ $\in$ $\Ac_{t,i}$. Thus, for any $\alpha\in\Ac^{K}_{t,i}(x)$, we have 
\beqs
 \E\big|N(\alpha)\big|^2 \; \leq \; \eta K(1+|x|^2), 
\enqs
which means that in the value functions $v_i(t,x)$ of optimal switching problems, one can restrict to controls $\alpha$ for which the second moment of $N(\alpha)$ is bounded by a constant depending on $x$.
}
\end{Remark}

\vspace{2mm}

Before proving Proposition \ref{proprestalpha}, we need the following Lemmata.

 \begin{Lemma}\label{unifboundXalpha}
For all $p$ $\geq$ $1$, there exists a positive constant $K_p$ such that
 \beqs
\sup_{\alpha\in\Ac_{t,i}} \Big\| \sup_{s\in[t,T]}\big|X_s^{t,x,\alpha}\big|\Big\|_{p} & \leq & K_p (1+|x|)\;, 
\enqs
for all $(t,x,i)$ $\in$ $[0,T]\times\R^d\times\I_q$. 
 \end{Lemma} 
 \textbf{Proof.} Fix $p$ $\geq$ $1$.  Then, we have from the definition of $X_s^{t,x,\alpha}$ in\reff{Xxalpha}, for  $(t,x,i)$ $\in$ $[0,T]\times\R^d\times\I_q$, $\alpha$ $\in$ $\Ac_{t,i}$: 
 \beqs
 \E\Big[ \sup_{s\in[t,r]}\big|X_s^{t,x,\alpha}\big|^p\Big]  & \leq & K_p\Big( |x|^p + \E\Big[\sum_{\tau_n\leq r} \int_{\tau_n}^{\tau_{n+1}\wedge r}\big|b_{\iota_n}(X_u^{t,x,\alpha})\big|^pdu\Big]\\
  & & + \; \E\Big[\sup_{s\in[t,r]}\Big|\sum_{\tau_n\leq s} \int_{\tau_n}^{\tau_{n+1}\wedge s}\sigma_{\iota_n}(X_u^{t,x,\alpha})dW_u\Big|^p\Big]\Big)\;, 
  \enqs
for all $r$ $\in$ $[t,T]$.  From the linear growth conditions on  $b_i$ and $\sigma_i$, for $i$ $\in$ $\I_q$, and Burkholder-Davis-Gundy's (BDG) inequality, we then get  by H\"older inequality when $p$ $\geq$ $2$:
 \beqs
 \E\Big[ \sup_{s\in[t,r]}\big|X_s^{t,x,\alpha}\big|^p\Big]  & \leq &  K_p\Big( 1+ |x|^p + \int_t^r\E\Big[\sup_{s\in[t,u]}\big|X_s^{t,x,\alpha}\big|^pdu\Big]\Big)\;,
  \enqs 
  for all $r$ $\in$ $[t,T]$. By  applying Gronwall's Lemma, we obtain the required estimate for $p$ $\geq$ $2$ , and then also for $p$ $\geq$ $1$ by H\"older inequality. 
 \ep

 \vspace{2mm}
 
\begin{Lemma}\label{lemcroislinv}
Under  {\bf (Hl)} and {\bf (Hc)}, the functions $v_i$, $i$ $\in$ $\I_q$,  satisfy a linear growth condition, i.e. there exists a constant $K$ such that 
\beqs
|v_i(t,x)| & \leq & K \big(1+|x|\big)\;,
\enqs
for all $(t,x,i)$ $\in$ $[0,T]\times\R^d\times\I_q$.
\end{Lemma}
\textbf{Proof.}  Under the linear growth condition on $f_i$, $g_i$ in {\bf (Hl)}, and the nonnegativity of the switching costs in {\bf (Hc)}, there exists some positive constant $K$ s.t.
\beqs
& & \E \Big[ \int_t^T f(X_s^{t,x,\alpha},I_s) ds + g(X_T^{t,x,\alpha},I_T)  - \sum_{n= 1}^{N(\alpha)} c(X_{\tau_n}^{t,x,\alpha},\iota_{n-1},\iota_n) \Big] \\ 
& \leq &  K \Big(1+\E\Big[\sup_{u\in[0,T]}\big|X^{t,x,\alpha}_u\big|\Big]\Big), 
\enqs
for all  $(t,x,i)$ $\in$ $[0,T]\times\R^d\times\I_q$, $\alpha$ $\in$ $\Ac_t,i$.  By combining with the estimate in Lemma  \ref{unifboundXalpha}, this shows that 
\beqs
v_i(t,x) & \leq & K(1+|x|)\;.
\enqs
Moreover,  by considering  the strategy $\alpha^0$ with no intervention i.e. $N(\alpha^0)$ $=$ $0$, we have
\beqs
v_i(t,x) & \geq & \E \Big[ \int_t^T f(X_s^{t,x,\alpha^0},i) ds + g(X_T^{t,x,\alpha^0},i) \Big]\\
 & \geq & - K \Big( 1+ \E\Big[\sup_{u\in[0,T]}\big|X^{t,x,\alpha}_u\big|\Big]\Big).
\enqs
Again, by the estimate in  Lemma \ref{unifboundXalpha}, this proves that 
\beqs
v_i(t,x) & \geq & -K(1+|x|)\;,
\enqs
and therefore the required linear growth condition on $v_i$.
\ep

\vspace{2mm}

\ni We now turn to the proof of  the Proposition.

 \vspace{2mm}
 
\ni \textbf{Proof of Proposition \ref{proprestalpha}.} The proof is done in 4 steps. Given $\alpha$ $\in$ $\Ac_{t,i}$, we will denote
\beqs
J(t,x,i;\alpha) &=& \E \Big[ \int_t^T f(X_s^{t,x,\alpha},I_s) ds + g(X_T^{t,x,\alpha},I_T) 
- \sum_{n= 1}^{N(\alpha)} c(X_{\tau_n}^{t,x,\alpha},\iota_{n-1},\iota_n) \Big].
\enqs

\ni $\bullet$ {\it Step 1}. First, we notice that the supremum in the definition of $v_i(t,x)$ may be taken over $\Ac_{t,i}^s$, where
\beqs
\Ac_{t,i}^s &=& \Big\{ \alpha = (\tau_n,\iota_n) \in \Ac_{t,i}:  (\tau_n) \mbox{ is strictly increasing }\Big\}.
\enqs
Indeed, it is always suboptimal to switch several times at a single date due to the triangular condition {\bf(Hc)}. 

\ni $\bullet$ {\it Step 2}. We now prove that it is enough to take the supremum over the strategies in $\Ac_{t,i}^{s,\infty}$, where
\beqs
\Ac_{t,i}^{s,\infty} &=& \Big\{ \alpha \in \Ac_{t,i}^s:  \E\big|C^{t,x,\alpha}_T\big|^2  \; < \; + \infty \Big\}.
\enqs
For any $\alpha = (\tau_k,\iota_k)_{k\geq 0}$ $\in$ $\Ac_{t,i}^s$, define $\alpha^n=(\tau_k^n,\iota_k^n)_{k\geq 0}$ as the strategy obtained from $\alpha$ by only keeping the first $n$ switches, i.e.
\beqs
 (\tau_k^n,\iota_k^n) &=& (\tau_k,\iota_k), \;\;\;\; k \leq n, \\
 \tau_k^n &=&  \infty, \;\;\;\;\;\;\;  k>n
\enqs
Note that for each $n$, $\alpha^n$ $\in$ $\Ac_{t,i}^{s,\infty}$.
Now since $\alpha$ and $\alpha^n$ (and the associated processes) coincide on $\left\{ N(\alpha) \leq n\right\}$, and by positivity of the switching costs,
\beqs
&&J(t,x,i;\alpha) - J(t,x,i;\alpha^n) \\ 
&\leq& \E \Big[ \Big(\int_t^T (f(X_s^{t,x,\alpha},I_s)-f(X_s^{t,x,\alpha^n},I_s)) ds + g(X_T^{t,x,\alpha},I_T) - g(X_T^{t,x,\alpha^n},I_T)
\Big) \mathbf{1}_{\{N(\alpha)>n\}} \Big] \\
&\leq& K (1+|x|) \P\big(N(\alpha)>n\big)^{1/2},
\enqs
by Cauchy-Schwarz inequality, linear growth of $f,g$ and Lemma \ref{unifboundXalpha}.
Hence letting $n$ $\rightarrow$ $\infty$, and since $N(\alpha)<\infty$ a.s., we obtain
\beqs
J(t,x,i;\alpha) \leq \liminf_{n \rightarrow \infty} J(t,x,i;\alpha^n),
\enqs
which proves the required assertion.

\ni$\bullet$ {\it Step 3}.  To each $\alpha$ $\in$ $\Ac_{t,i}^{s,\infty}$, we  associate the process $(Y^{t,x,\alpha},Z^{t,x,\alpha})$ solution to the following Backward Stochastic Differential Equation (BSDE)
\beq\label{EDSRalpha}
Y^{t,x,\alpha}_u & = & g(X_T^{t,x,\alpha},I^\alpha_T)+ \int_u^T f(X_s^{t,x,\alpha},I^\alpha_s) ds \\
& & \;\;\; - \; \int_u^T Z^{t,x,\alpha}_sdW_s - C^{t,x,\alpha}_T + C^{t,x,\alpha}_u\;,\qquad t\leq u\leq T\nonumber
\enq
and satisfying the condition
\beqs
\E\Big[ \sup_{s\in[t,T]}  |Y^{t,x,\alpha}_s|^2 \Big]+\E\Big[ \int_{t}^T  |Z_s^{t,x,\alpha}|^2ds \Big] & < & \infty.  
\enqs
Such a solution exists under {\bf (Hl)},  Lemma \ref{unifboundXalpha} and $\E\big[|C^{t,x,\alpha}_T|^2\big]$ $<$ $\infty$.
Note that taking the expectation in \reff{EDSRalpha}, $Y^{t,x,\alpha}_t$ $=$ $J(t,x,i;\alpha)$.

We now define for $\tilde K>0$,
\beqs
\tilde{\Ac}_{t,i}^{s,\tilde K}(x) &=& \Big\{ \alpha \in \Ac_{t,i}^{s,\infty}:  \E\Big[\sup_{s \in [t,T]}\big|Y^{t,x,\alpha}_s\big|^2\Big]  \; \leq \; \tilde{K}(1+|x|^2) \Big\},
\enqs
and claim that for some constant $\tilde K$, the supremum in $v_i(t,x)$ may be taken over $\alpha$ $\in$ $\tilde{\Ac}_{t,i}^{s,\tilde K}(x)$.
First taking the conditional expectation in \reff{EDSRalpha}, we have
\beqs
Y^{t,x,\alpha}_u \leq v_{I_t}(X_u^{t,x,\alpha},I^\alpha_u) \leq K(1+\big|X_s^{t,x,\alpha}\big|), \;\;\; t \leq u \leq T, 
\enqs
so that by Lemma \ref{unifboundXalpha} the only restriction is to have a lower bound on $Y^{t,x,\alpha}_u$. As in Lemma \ref{lemcroislinv}, this is done by considering strategies with fewer interventions.
%
Given $\alpha$ $\in$ $\Ac_{t,i}^{s,\infty}$, consider the stopping time 
\beqs
\tau = \inf \{ s\geq t :  J(s,X_s^{t,x,\alpha},I_s^\alpha;\alpha^0) \geq Y^{t,x,\alpha}_s \}
\enqs
where $\alpha^0$ is the strategy with no switches, and define $\tilde{\alpha} = (\tilde{\tau}_n, \iota_n)$, where
\beqs
\tilde{\tau}_n &=& \tau_n \mathbf{1}_{\{\tau_n \leq \tau\}} + \infty \mathbf{1}_{\{\tau_n > \tau\}}.
\enqs
Now for each $t \leq u \leq T$, taking the conditional expectation in \reff{EDSRalpha} we obtain
\beqs
&&\mathbf{1}_{\{u\leq\tau\}}(Y_u^{t,x,\tilde \alpha} - Y_u^{t,x,\alpha})  \\
&=&  \E \Big[\mathbf{1}_{\{u\leq\tau<T\}} \Big(\int_\tau^T f(X_s^{t,x,\tilde \alpha},I_s^{\tilde\alpha}) ds  + g(X_T^{t,x,\tilde \alpha},I_T^{\tilde\alpha}) \\
&& \;\;\;\;\;\;  - \int_\tau^T f(X_s^{t,x,\alpha},I_s) ds  - g(X_T^{t,x,\alpha},I_T) + C^{t,x,\alpha}_T - C^{t,x,\alpha}_\tau \Big) \big| \Fc_u \Big] \\
&=& \E \Big[\mathbf{1}_{\{u\leq\tau<T\}} \big( J(\tau, X_{\tau}^{t,x,\alpha}, I_\tau^{\alpha};\alpha^0) - Y_\tau^{t,x,\alpha}\big) \big |\Fc_u \Big],
\enqs
where we have taken the conditional expectation w.r.t. $\Fc_\tau$ inside the expectation. Since the process $\big(J(u, X_{u}^{t,x,\alpha}, I_u^{\alpha};\alpha^0) - Y_u^{t,x,\alpha}\big)_{t\leq u \leq T}$ has right-continuous paths , by definition of $\tau$ we have  $J(\tau, X_{\tau}^{t,x,\alpha}, I_\tau^{\alpha},\alpha^0) - Y_\tau^{t,x,\alpha}$ $\geq$ $0$ a.s., so that
\beq \label{diffY-tildeY}
\mathbf{1}_{\{u\leq\tau\}}(Y_u^{t,x,\tilde \alpha} - Y_u^{t,x,\alpha})  & \geq &  0\;.
\enq
Noting that on $\{u\leq \tau\}$ we have
\beqs
Y_u^{t,x,\alpha} &=&  Y_{u-}^{t,x,\alpha} + \Delta Y_u^{t,x,\alpha} \\
&\geq& J(u,X_u^{t,x,\alpha},I_{u-}^\alpha; \alpha^0_u) +  c(X_u^{t,x,\alpha},I_{u-}^\alpha,I_u^\alpha) \\
&\geq& - K (1+ |X_u|)\;,
\enqs
and since on $\{u > \tau\}$, $Y_u^{t,x,\tilde \alpha} = J(u,X^{t,x,\tilde \alpha}_u,I_u^{\tilde\alpha};\alpha^0)$, from Lemma \ref{unifboundXalpha},  it follows that $\tilde \alpha$ $\in$ $\tilde \Ac^{s,\tilde K}_{t,i}(x)$, for some 
$\tilde K$ not depending on $(t,x)$. Furthermore taking $u=t$ in \reff{diffY-tildeY}, we have $J(t,x,i;\tilde \alpha) \geq J(t,x,i; \alpha)$, and this  proves the required assertion.

\ni$\bullet$ {\it Step 4}.  Finally we show that for each $\tilde K$, there exists some positive $K$ s.t. $\tilde{\Ac}_{t,i}^{s,\tilde K}(x)$ $\subset$ $\Ac_{t,i}^{K}(x)$. 
We fix $\alpha$ $\in$ $\tilde{\Ac}_{t,i}^{s,\tilde K}(x)$. Applying It\^o's formula to $|Y^{t,x,\alpha}|^2$ in \reff{EDSRalpha}, we have 
\beqs
|Y^{t,x,\alpha}_t|^2+\int_t^T|Z^{t,x,\alpha}_s|^2ds & = & |g(X_T^{t,x,\alpha},I^\alpha_T)|^2+2 \int_t^T Y^{t,x,\alpha}_sf(X_s^{t,x,\alpha},I^\alpha_s) ds\\
 & &\;  - \; 2\int_t^TY^{t,x,\alpha}_sZ^{t,x,\alpha}_sdW_s- 2\int_t^T Y^{t,x,\alpha}_s dC^{t,x,\alpha}_s\;.
\enqs
Using \textbf{(Hl)} and the inequality $2ab$ $\leq$ $a^2+b^2$ for $a,b$ $\in$ $\R$, we get 
\beq\nonumber
\int_t^T|Z^{t,x,\alpha}_s|^2ds  & \leq  &  K \Big(1+\sup_{s\in[t,T]}|X_s^{t,x,\alpha}|^2+  \sup_{s\in[t,T]}|Y_s^{t,x,\alpha}|^2 +  |C^{t,x,\alpha}_T-C^{t,x,\alpha}_t|\sup_{s\in[t,T]}|Y_s^{t,x,\alpha}| \Big) \\
 & &  \;\;\;\;\;\;\;  - 2\int_t^TY^{t,x,\alpha}_sZ^{t,x,\alpha}_sdW_s\;.\label{est1}
\enq
Moreover, from \reff{EDSRalpha}, we have
\beq\nonumber
|C^{t,x,\alpha}_T-C^{t,x,\alpha}_t|^2 & \leq & K\Big( 1+\sup_{s\in[t,T]}|X_s^{t,x,\alpha}|^2+\sup_{s\in[t,T]}|Y_s^{t,x,\alpha}|^2 \\
 & & \;\;\; \;\;\;\; + \; \Big|\int_t^TZ^{t,x,\alpha}_sdW_s\Big|^2 \Big)\label{est2}
\enq 
Combining \reff{est1} and \reff{est2} and using the inequality $ab$ $\leq$ ${a^2\over 2\eps}+{\eps b^2\over 2}$, for $a,b$ $\in$ $\R$ and $\eps$ $>$ $0$, we obtain 
\beqs
\int_t^T|Z^{t,x,\alpha}_s|^2ds &\leq &  K\Big((1+\eps)\Big(1+\sup_{s\in[t,T]}|X_s^{t,x,\alpha}|^2\Big)+\sup_{s\in[t,T]}|Y_s^{t,x,\alpha}|^2\big(\eps+{1\over\eps}\big)\\
 & & \;\;\;\;\;\;\;    + \; \eps \Big|\int_t^TZ^{t,x,\alpha}_sdW_s\Big|^2 \Big) - 2 \int_t^TY^{t,x,\alpha}_sZ^{t,x,\alpha}_sdW_s\;.
\enqs
Taking the expectation in the previous estimate, it follows from Lemma \ref{unifboundXalpha} and $\alpha$ $\in$ $\tilde{\Ac}_{t,i}^{s,\tilde K}(x)$ that
\beq
\E\Big[\int_t^T|Z^{t,x,\alpha}_s|^2ds\Big] & \leq & K \Big((1+\eps)\Big(1+\E\sup_{s\in[t,T]}|X_s^{t,x,\alpha}|^2\Big)+\big(\eps+{1\over\eps}\big)\E\sup_{s\in[t,T]}|Y_s^{t,x,\alpha}|^2 \nonumber \\
 & &  \; + \; \eps\E \Big|\int_t^TZ^{t,x,\alpha}_sdW_s\Big| \nonumber \\
  & \leq &  K \Big((1+|x|^2)\big(1+\eps+{1\over\eps}\big)+\eps \E\Big[\Big(\int_t^T|Z^{t,x,\alpha}_s|^2ds\Big)\Big]\Big)\;,  \nonumber 
\enq
Taking $\eps$ small enough, this yields
\beqs
\E\Big[\int_t^T|Z^{t,x,\alpha}_s|^2ds\Big] & \leq &  K \big(1+|x|^2\big)\;,
\enqs 
Taking the expectation in \reff{est2}, and  using the previous inequality together with Lemma \ref{unifboundXalpha} and $\alpha$ $\in$ $\tilde{\Ac}_{t,i}^{s,\tilde K}(x)$, we get: 
\beq \label{estimC}
\E |C^{t,x,\alpha^*}_T -C^{t,x,\alpha^*}_t|^2 & \leq & K(1+|x|^2),
\enq
for some positive constant $K$ not depending on $(t,x,i)$. 
Since $(\tau_n)$ is strictly increasing, we know that at the initial time $t$, there is at most one decision time $\tau_1$.  Thus, from the linear growth condition on the 
switching cost, $\E[|C_t^{t,x,\alpha}|^2]\leq K(1+|x|^2)$, which implies with \reff{estimC} that $\alpha$ $\in$ $\Ac_{t,i}^{K}(x)$, and this proves the required result.
\ep

\vspace{3mm}

 In the sequel of this paper, we shall assume that {\bf (Hl)} and {\bf (Hc)} stand in force.

\section{Time discretization}

\setcounter{equation}{0} \setcounter{Assumption}{0}
\setcounter{Theorem}{0} \setcounter{Proposition}{0}
\setcounter{Corollary}{0} \setcounter{Lemma}{0}
\setcounter{Definition}{0} \setcounter{Remark}{0}

We first consider a time discretization of $[0,T]$ with time step $h$ $=$ $T/m$ $\leq$ $1$, and partition $\T_h$ $=$ $\{t_k=kh, k=0,\ldots,m\}$. 
 For $(t_k,i)$ $\in$ $\T_h\times\I_q$, we denote by $\Ac_{t_k,i}^h$ the set of admissible switching controls $\alpha$ $=$ $(\tau_n,\iota_n)_n$ in 
 $\Ac_{t_k,i}$, such that $\tau_n$ are  valued in $\{\ell h, \ell = k,\ldots,m\}$, and we consider the value functions for the discretized optimal switching problem:
 \beq 
v_i^h(t_k,x) &=& \sup_{\alpha\in\Ac_{t_k,i}^h}   \E \Big[ \sum_{\ell=k}^{m-1} f(X_{t_\ell}^{t_k,x,\alpha},I_{t_\ell}) h  + g(X_{t_m}^{t_k,x,\alpha},I_{t_m}) 
\nonumber \\ 
& & \hspace{3cm} - \;  \sum_{n=1}^{N(\alpha)} c(X_{\tau_n}^{t_k,x,\alpha},\iota_{n-1},\iota_n) \Big], \label{defvih}
\enq
for $(t_k,i,x)$ $\in$ $\T_h\times\I_q\times\R^d$.   
 
\vspace{2mm}

The next result  provides an error analysis between the continuous-time optimal switching problem and its discrete-time version.

\begin{Theorem}\label{diff-v-v^h} 
There exists a positive constant $K$ (not depending on $h$) such that 
\beqs
|v_i(t_k,x) - v_i^h(t_k,x) | & \leq & K(1 + |x|^{5/2}) \left(h \log(2T/h)\right)^{1/2}, 
\enqs 
for  all $(t_k,x,i)$  $\in$ $\T_h\times \R^d\times\I_q$. 

If the cost functions $c_{ij}$, $i,j$ $\in$ $\I_q$, do not depend on $x$, then
\beqs
|v_i(t_k,x) - v_i^h(t_k,x) | & \leq & K(1 + |x|^{3/2}) h^{1/2}
\enqs
\end{Theorem}

\vspace{1mm}

\begin{Remark} 
{\rm  For optimal stopping problems, it is known that the approximation by the  discrete-time version gives an 
error of order  $h^{1\over 2}$, see e.g. \cite{lam02} and \cite{balpag03}.   We recover this rate of convergence for multiple switching problems when the switching costs do not depend on the state process. However, in the general case, the error is of order $(h \log(1/h))^{{1\over 2}}$. A rate of $h^{{1 \over 2} - \varepsilon}$ was obtained in \cite{chaelikha10} in the case of uncontrolled state process $X$, and is improved and extended here when $X$ may be influenced  through its drift and diffusion coefficient by the switching control. 
}
\end{Remark}

\vspace{1mm}

Before proving this Theorem, we need the three following lemmata. The first two deal with the regularity in time of the controlled diffusion uniformly in the control, and the third one deals 
with the regularity of the controlled diffusion with respect to  the control. 

\begin{Lemma}\label{lem-reg-unif}
There exists a constant $K$ such that 
\beqs
\sup_{\alpha\in\Ac_{t_k,i}}  \max_{k\leq \ell\leq m-1}
 \Big\|\sup_{s\in[t_\ell,t_{\ell+1}]}  \big|X_s^{t_k,x,\alpha}-X_{t_\ell}^{t_k,x,\alpha}\big| \Big\|_2    & \leq & K(1+ |x|){h}^{1\over 2}\;,
\enqs
for all $x$ $\in$ $\R^d$, $i$ $\in$ $\I_q$, $k$ $=$ $0,\ldots,n$.
\end{Lemma}
\textbf{Proof.} From the definition of $X^{t,x,\alpha}$ in \reff{Xxalpha}, we have for all 
$(t_k,x,i)$ $\in$ $\T_h\times\R^d\times\I_q$ and $\alpha$ $\in$ $\Ac_{t_k,i}$, 
 \beqs
 \E\Big[\sup_{u\in[t_\ell, s]} \big|X_u^{t,x,\alpha}-X_{t_\ell}^{t,x,\alpha}\big|^2\Big]  & \leq & 
 K\Big(\E\Big[\Big(\int_{t_\ell}^{s} |b_{I_u}(X_u^{t,x,\alpha})|du\Big)^2\Big]  \\
 & & \;\;\;\;\;\;\;\;  + \;  
 \E\Big[\sup_{u\in[t_\ell, s]}\Big| \int_{t_\ell}^{ u}\sigma_{I_r}(X_r^{t,x,\alpha})dW_r\Big|^2\Big]\Big)\;,
  \enqs
  for all $s$ $\in$ $[t_\ell,t_{\ell+1}]$. From  BDG  and Jensen inequalities, we then have
  \beqs
 \E\Big[\sup_{u\in[t_\ell, s]} \big|X_u^{t,x,\alpha}-X_{t_\ell}^{t,x,\alpha}\big|^2\Big]  & \leq & K \Big(\E\Big[\int_{t_\ell}^{s}\big|b_{I_u}(X_u^{t,x,\alpha})\big|^2du\Big]
+\E\Big[\int_{t_\ell}^{ s}\big| \sigma_{I_u}(X_u^{t,x,\alpha})\big|^2du\Big]\Big)\;,
  \enqs
From the linear growth conditions on $b_i$ and $\sigma_i$, for $i$ $\in$ $\I_q$, and Lemma \ref{unifboundXalpha}, we conclude that
 \beqs
\E\Big[\sup_{s\in[t_\ell,t_{\ell+1}]}\big|X_s^{t,x,\alpha}-X_{t_\ell}^{t,x,\alpha}\big|^p\Big]  & \leq &  K_p( 1+ |x|^p)h.
\enqs
 \ep

\begin{Lemma} \label{lem-reg-unif-2}
There exists some positive  constant $K$ such that 
\beqs
\sup_{\alpha\in\Ac_{t_k,i}}  \Big\|\sup_{\substack{0\leq s,u\leq T \\ |s-u| \leq h}}  \big|X_s^{t_k,x,\alpha}-X_{u}^{t_k,x,\alpha}\big| \Big\|_2    & \leq & K(1+ |x|)\big(h \log(2T/h) \big)^{1\over 2}\;,
\enqs
\end{Lemma}
\textbf{Proof.} This follows from Theorem 1 in \cite{fisnap10}, using the estimates from Lemma \ref{unifboundXalpha} and linear growth of $b_i$, $\sigma_i$.
\ep

\vspace{2mm}

For a strategy $\alpha = (\tau_n,\iota_n)_n$ $\in$ $\Ac_{t_k,i}$ we denote by $\tilde \alpha = (\tilde  \tau_n,\tilde \iota_n)_n$ the strategy of  $\Ac_{t_k,i}^h$  defined by
\beqs
\tilde  \tau_n \; = \;  \min\{t_\ell \in \T_h~:~t_\ell\geq \tau_n\}\;, & & \tilde \iota_n \; = \; \iota_n, \;\;\;\;\;\; n \in \N. 
\enqs 
The strategy $\tilde \alpha$ can be seen as the approximation of the strategy $\alpha$ by an element of $\Ac^h_{t_k,i}$. We then have the following regularity result of the diffusion in the control $\alpha$.

\begin{Lemma}\label{lem-diff-x-xtilde}
There exists some positive  constant $K$ such that 
\beqs
 \Big\| \sup_{s\in[t_k,T]}  \big|X_s^{t_k,x,\alpha}-X_{s}^{t_k,x,\tilde \alpha}\big| \Big\|_{2}   & \leq & K \Big(\E[N(\alpha)^2]\Big)^{1\over 4} (1+ |x|){h}^{1\over 2}\;,
\enqs
for all $x$ $\in$ $\R^d$, $i$ $\in$ $\I_q$, $k$ $=$ $0,\ldots,n$ and $\alpha$ $\in$ $\Ac_{t_k,i}$. 
\end{Lemma}
\textbf{Proof. }
From the definition of $X^{t,x,\alpha}$ and $X^{t,x,\tilde \alpha}$, for  $(t_k,x,i)$ $\in$ $\T_h\times\R_d\times\I_q$, $\alpha$ $\in$ $\Ac_{t_k,i}^K$, we have by BDG inequality:  
\beqs
\E\Big[\sup_{u\in[t_k,s]}\big|X_u^{t,x,\alpha}-X_{u}^{t,x,\tilde \alpha}\big| ^2  \Big] & \leq & K \Big( \E\Big[\int_{t_k}^{s}\big|b_{}(X_u^{t,x,\alpha},I_u)-b_{}(X_u^{t,x,\tilde \alpha},\tilde I_u)\big|^2du\Big]\\
 & & \;\; + \; \E\Big[ \int_{t_k}^{ s}\big|\sigma_{}(X_u^{t,x,\alpha},I_u)-\sigma(X_u^{t,x,\tilde \alpha},\tilde I_u)\big|^2du\Big] \Big)\;,
\enqs
for all $s$ $\in$ $[t_k,T]$. Then using Lipschitz property of $b_i$ and $\sigma_i$ for $i$ $\in$ $\I_q$ we get: 
\beq \nonumber
\E\Big[\sup_{u\in[t_k,s]}\big|X_s^{t,x,\alpha}-X_{s}^{t,x,\tilde \alpha}\big| ^2  \Big] & \leq & K  \Big( \E\Big[\int_{t_k}^{s}\big|X_u^{t,x,\alpha}-X_u^{t,x,\tilde \alpha}\big|^2du\Big]\\\nonumber
 & & \;\; 
+ \; \E\Big[ \int_{t_k}^{ s}\big|b(X_u^{t,x,\alpha}, I_u)-b(X_u^{t,x, \alpha},\tilde I_u)\big|^2du\Big] \\\nonumber
 & &  \;\; + \; \E\Big[ \int_{t_k}^{ s}\big|\sigma_{}(X_u^{t,x,\alpha},I_u)-\sigma(X_u^{t,x, \alpha},\tilde I_u)\big|^2du\Big] \Big)\\\label{ineqprec}
& \leq & K\Big( \E\Big[\int_{t_k}^{s}\sup_{r\in[t_k,u]}\big|X_r^{t,x,\alpha}-X_r^{t,x,\tilde \alpha}\big|^2du\Big]\\\nonumber
 & & \;\;\;  + \; \E\Big[\big(\sup_{u\in[t_k,T]}\big|X_u^{t,x,\alpha}\big|^2+1\big)\int_{t_k}^{s}\mathbf{1}_{I_s\neq \tilde I_s} ds\Big]\Big)\;,
\enq
for all $s$ $\in$ $[t_k,T]$.  From the definition of $\tilde \alpha$ we have 
\beqs
\int_{t_k}^{s}\mathbf{1}_{I_s\neq \tilde I_s}ds & \leq &  N(\alpha) h\;,
\enqs
which gives with \reff{ineqprec}, Lemma \ref{unifboundXalpha}, Remark \ref{remNal} and H\"older inequality: 
\beqs
\E\Big[\sup_{u\in[t_k,s]}\big|X_u^{t,x,\alpha}-X_{u}^{t,x,\tilde \alpha}\big| ^2  \Big] & \leq &
K\Big( \E\Big[\int_{t_k}^{s}\sup_{r\in[t_k,u]}\big|X_r^{t,x,\alpha}-X_r^{t,x,\tilde \alpha}\big|^2du\Big] \\
& & \;\;\;\;\;\;\;\;  + \;  \big(\E[N(\alpha)^2]\big)^{1\over 2} (1+ |x|^2)h\Big),
\enqs
for all $s$ $\in$ $[t_k,T]$.  We conclude with   Gronwall's Lemma.
\ep

\vspace{2mm}

\ni We are now ready to prove the convergence result for the time discretization of the optimal switching problem.

\vspace{2mm}

\ni\textbf{Proof of Theorem \ref{diff-v-v^h}.} We introduce the auxiliary function $\tilde v^h_i$ defined by
\beqs
\tilde v ^h_i(t_k,x) & = & \sup_{\alpha\in\Ac_{t_k,i}^h}   \E \Big[ \int_{t_k}^T f(X_s^{t_k,x,\alpha},I_s) ds + g(X_T^{t_k,x,\alpha},I_T) 
- \sum_{n= 1}^{N(\alpha)} c(X_{\tau_n}^{t_k,x,\alpha},\iota_{n-1},\iota_n) \Big]\;,
\enqs
for all $(t_k,x)$ $\in$ $\T_h\times\R^d$. We then write 
\beqs
| v_i(t_k,x) - v_i^h(t_k,x) | & \leq & | v_i(t_k,x) - \tilde v_i^h(t_k,x) |  + |  \tilde v_i^h(t_k,x) - v_i^h(t_k,x) | \;,
\enqs
and study each of the two terms in the right-hand side. 

\ni $\bullet$ Let us investigate the  first term.  By definition  of  the approximating strategy $\tilde\alpha$ $=$ $(\tilde\tau_n,\tilde\iota_n)_n$ $\in$ 
$\Ac_{t_k,i}^h$ of  $\alpha$ $\in$ $\Ac_{t_k,i}$, we see that the auxiliary value function $\tilde v_i^h$ may be written as
\beqs
\tilde v ^h_i(t_k,x) & = & \sup_{\alpha\in\Ac_{t_k,i}}   \E \Big[ \int_{t_k}^T f(X_s^{t_k,x,\tilde \alpha},\tilde I_s) ds + g(X_T^{t_k,x,\tilde \alpha},\tilde I_T) 
- \sum_{n= 1}^{N(\alpha)} c(X_{\tilde \tau_n}^{t_k,x,\tilde \alpha},\tilde \iota_{n-1},\tilde \iota_n) \Big],
\enqs
where $\tilde I$ is the indicator of the regime value associated to $\tilde\alpha$. 
 Fix now a positive number $\bar K$ s.t. relation \reff{relviK}  in Proposition \ref{proprestalpha} holds, and observe that
 \beqs
& &  \sup_{\alpha\in\Ac^{\bar K}_{t_k,i}(x)}   \E \Big[ \int_{t_k}^T f(X_s^{t_k,x,\tilde \alpha},\tilde I_s) ds + g(X_T^{t_k,x,\tilde \alpha},\tilde I_T) 
- \sum_{n= 1}^{N(\alpha)} c(X_{\tilde \tau_n}^{t_k,x,\tilde \alpha},\tilde \iota_{n-1},\tilde \iota_n) \Big]  \\
& \leq & \tilde v ^h_i(t_k,x) \;  \leq \;  v_i(t_k,x) \\
&=&  \sup_{\alpha\in\Ac^{\bar K}_{t_k,i}(x)}   \E \Big[ \int_{t_k}^T f(X_s^{t_k,x,\alpha}, I_s) ds + g(X_T^{t_k,x,\alpha},I_T) 
- \sum_{n= 1}^{N(\alpha)} c(X_{\tau_n}^{t_k,x,\alpha},\iota_{n-1},\iota_n) \Big].
\enqs
We then have
\beq \label{videlta}
|v_i(t_k,x) - \tilde v_i^h(t_k,x) | & \leq & \sup_{\alpha\in\Ac^{\bar K}_{t_k,i}(x)}  \Big[ \Delta^1_{t_k,x}(\alpha)+\Delta^2_{t_k,x}(\alpha) \Big],
\enq
with 
\beqs
\Delta^1_{t_k,x}(\alpha) & = &  \E \Big[ \int_{t_k}^T \big|f(X_s^{t_k,x, \alpha}, I_s) -f(X_s^{t_k,x,\tilde \alpha},\tilde I_s)\big|ds + \big|g(X_T^{t_k,x, \alpha}, I_T) -g(X_T^{t,x,\tilde \alpha},\tilde I_T) \big|\Big] \;,\\
\Delta^2_{t_k,x}(\alpha) & = &\E\Big[\sum_{n= 1}^{N(\alpha)} \big|c(X_{ \tau_n}^{t_k,x, \alpha}, \iota_{n-1}, \iota_n)-c(X_{\tilde \tau_n}^{t_k,x,\tilde \alpha},\tilde \iota_{n-1},\tilde \iota_n)\big| \Big].
\enqs
Under {\bf (Hl)}, and by definition of $\tilde \alpha$,  there exists some positive constant $K$ s.t. 
\beq
\Delta^1_{t_k,x}(\alpha) & \leq & K \Big( \sup_{s\in[t_k,T]}\E \Big[  \big|X_s^{t_k,x, \alpha} -X_s^{t_k,x,\tilde \alpha}\big|\Big] + \E \Big[\big(\sup_{s\in[t_k,T]}\big|X_s^{t_k,x, \alpha}\big|+1\big)\int_{t_k}^T\mathbf{1}_{I_s\neq \tilde I_{s}}ds \Big] \Big) \;.~\qquad\nonumber \\ \label{majA0}
 & \leq &  K \Big( \sup_{s\in[t_k,T]}\E \Big[  \big|X_s^{t_k,x, \alpha} -X_s^{t_k,x,\tilde \alpha}\big|\Big] \\
  & & \;\;\;\;\; + \;  \Big(1+ \Big\| \sup_{s\in[t_k,T]}\big|X_s^{t_k,x, \alpha}\big| \Big\|_{2} \Big)
  \Big(\E \Big[\int_{t_k}^T\mathbf{1}_{I_s\neq \tilde I_{s}}ds \Big] \Big)^{1\over 2} \Big),\nonumber
\enq
by Cauchy-Schwarz inequality.  For  $\alpha$ $\in$ $\Ac^{K}_{t_k,i}(x)$, we have by Remark \ref{remNal}
\beqs
\E \Big[\int_{t_k}^T\mathbf{1}_{I_s\neq \tilde I_{s}}ds\Big] & \leq & h\E \Big[N(\alpha)\Big]  \;  \leq \;   
\eta \bar K_1 (1+|x|)h,
\enqs
for some positive constant $\eta$ $>$ $0$. By using this last estimate together with  Lemmata \ref{unifboundXalpha} and  \ref{lem-diff-x-xtilde} into \reff{majA0}, we obtain the existence of some constant $K$ s.t.
\beq\label{majA}
\sup_{\alpha \in \Ac_{t_k,i}^{\bar K}(x)}\Delta^1_{t_k,x}(\alpha) & \leq &  K(1+|x|^{3/2})h^{1\over 2},
\enq
for all  $(t_k,x,i)$ $\in$ $\T_h\times\R^d\times\I_q$. 

We now turn to the term $\Delta^2_{t,x}(\alpha)$. Under {\bf (Hl)}, and by definition of $\tilde \alpha$,  there exists some positive constant $K$ s.t.
\beq
\Delta^2_{t_k,x}(\alpha) & \leq &  K \E\Big[\sum_{n= 1}^{N(\alpha)} \big|X_{ \tau_n}^{t_k,x, \alpha}-X_{\tilde \tau_n}^{t_k,x,\tilde \alpha}\big| \Big]\nonumber \\
 & \leq & K \Big(\E\Big[\sum_{n= 1}^{N(\alpha)} \big|X_{ \tau_n}^{t_k,x, \alpha}-X_{\tilde \tau_n}^{t_k,x, \alpha}\big|\Big] +
 \E\Big[{N(\alpha)} \sup_{s\in[t_k,T]}\big|X_{s}^{t_k,x, \alpha}-X_{s}^{t_k,x,\tilde \alpha}\big|\Big]\Big)\nonumber\\
  & \leq & K \Big(\E\Big[\sum_{n= 1}^{N(\alpha)} \big|X_{ \tau_n}^{t_k,x, \alpha}-X_{\tilde \tau_n}^{t_k,x, \alpha}\big|\Big]  \nonumber \\
  & &\;\;\;\;\; \;   + \;  \Big\| N(\alpha) \Big\|_{2}   \Big\| \sup_{s\in[t_k,T]}\big|X_{s}^{t_k,x, \alpha}-X_{s}^{t_k,x,\tilde \alpha}\big| \Big\|_{2} \Big), 
   \label{decomp Delta2}
\enq
by Cauchy-Schwarz inequality. For $\alpha$ $\in$ $\Ac_{t_k,i}^{K}(x)$ with Remark  \ref{remNal},  and from Lemma \ref{lem-diff-x-xtilde}, 
we get the existence of some positive constant $K$ s.t. 
\beq\label{morceau1Delta2}
\Big\| N(\alpha) \Big\|_{2}   \Big\| \sup_{s\in[t_k,T]}\big|X_{s}^{t_k,x, \alpha}-X_{s}^{t_k,x,\tilde \alpha}\big| \Big\|_{2} & \leq & 
K(1+|x|^{5/2})h^{1\over 2}\;.
 \enq
On the other hand, 
\beqs
\E\Big[\sum_{n= 1}^{N(\alpha)} \big|X_{ \tau_n}^{t_k,x, \alpha}-X_{\tilde \tau_n}^{t_k,x, \alpha}\big|\Big] & \leq & 
\E\Big[ N(\alpha) \sup_{\substack{0\leq s,u\leq T \\ |s-u| \leq h}}  \big|X_s^{t_k,x,\alpha}-X_{u}^{t_k,x,\alpha}\big|\Big] \\
&\leq&  \Big\|N(\alpha)\Big\|_2 \Big\|\sup_{\substack{0\leq s,u\leq T \\ |s-u| \leq h}}  \big|X_s^{t_k,x,\alpha}-X_{u}^{t_k,x,\alpha}\big| \Big\|_2
\enqs
by Cauchy-Schwarz inequality.  For $\alpha$ $\in$ $\Ac^{\bar K}_{t_k,i}(x)$, by  Lemma \ref{lem-reg-unif-2}, this yields  the existence of some positive constant $K$ s.t. 
\beq\label{morceau2Delta2}
\E\Big[\sum_{n= 1}^{N(\alpha)} \big|X_{ \tau_n}^{t_k,x, \alpha}-X_{\tilde \tau_n}^{t_k,x, \alpha}\big|\Big] & \leq & 
K(1+|x|^2) \left(h \log(2T/h)\right)^{1/2}. 
\enq 
By plugging \reff{morceau1Delta2} and \reff{morceau2Delta2} into  \reff{decomp Delta2}, we then get 
\beq\label{majB}
\Delta^2_{t,x}(\alpha) & \leq & K(1+|x|^2) \left(h \log(2T/h)\right)^{1/2}\;.
\enq
Combining \reff{majA} and \reff{majB}, we obtain with \reff{videlta}
\beqs
| v_i(t_k,x) - \tilde v_i^h(t_k,x) | & \leq & K(1+|x|^2) \left(h \log(2T/h)\right)^{1/2}\;.
\enqs
In the case where $c$ does not depend on the variable $x$, we have $\Delta^2_{t,x}(\alpha)$ $=$ $0$, and so by  \reff{videlta}, \reff{majA}: 
\beqs
| v_i(t_k,x) - \tilde v_i^h(t_k,x) | & \leq & K(1+|x|^{3/2})h^{{1\over 2}}\;.
\enqs

\vspace{2mm}

\ni$\bullet $ For the second term, we have by definition of $v_i^h$ and $\tilde v_i^h$: 
\beqs
|  \tilde v_i^h(t_k,x) - v_i^h(t_k,x) | & \leq & \sup_{\alpha\in\Ac_{t_k,i}^h}   \E \Big[ \sum_{\ell=k}^{m-1} 
\int_{t_\ell}^{t_{\ell+1}} \big|f(X_s^{t,x,\alpha},I_s)-f(X_{t_\ell}^{t,x,\alpha},I_s)\big| ds   \Big],
\enqs
since $I_s$ $=$ $I_{t_\ell}$ on $[t_\ell,t_{\ell+1})$. 
Under \textbf{(Hl)}, we get
\beqs
|  \tilde v_i^h(t_k,x) - v_i^h(t_k,x) | & \leq &K \sup_{\alpha\in\Ac_{t_k,i}^h} \max_{k\leq \ell \leq m-1} \sup_{s\in[t_\ell,t_{\ell+1}]}  
 \E \Big[\big|X_s^{t,x,\alpha}-X_{t_\ell}^{t,x,\alpha}\big|    \Big],
\enqs
for some positive constant $K$, and by Lemma \ref{lem-reg-unif}, this shows that
\beqs
| \tilde v_i^h(t_k,x) - v_i^h(t_k,x) | & \leq & K(1+|x|){h}^{1\over 2}\;.
\enqs
\ep

\vspace{2mm}

In  a second  step, we approximate the continuous-time (controlled) diffusion by a discrete-time (controlled) Markov chain  following an Euler type scheme.  For any $(t_k,x,i)$ $\in$ $\T_h\times\R^d\times\I_q$, $\alpha$ $\in$ $\Ac_{t_k,i}^h$, we introduce  
$(\bar X_{t_\ell}^{h,t_k,x,\alpha})_{k\leq \ell\leq m}$ defined by: 
\beqs
\bar X_{t_k}^{h,t_k,x,\alpha} \; =\;  x, & & \bar X_{t_{\ell+1}}^{h,t_k,x,\alpha} \; = \;  F_{I_{t_\ell}}^h(\bar X_{t_\ell}^{h,t_k,x,\alpha},\vartheta_{\ell+1}), 
\;\;\;\; k \leq \ell \leq m-1,
\enqs 
where 
\beqs
 F_i^h(x,\vartheta_{k+1}) \; = \;  x + b_i(x) h + \sigma_i(x) \sqrt{h} \; \vartheta_{k+1},
\enqs
and $\vartheta_{k+1}$ $=$  $(W_{t_{k+1}}-W_{t_k})/\sqrt{h}$, $k$ $=$ $0,\ldots,m-1$,  are iid, $\Nc(0,I_d)$-distributed, independent of $\Fc_{t_k}$.  
Similarly as in Lemma \ref{unifboundXalpha}, we have the $L^p$-estimate: 
\beq \label{estimEuler}
\sup_{\alpha\in\Ac^h_{t_k,i}} \Big\| \max_{\ell = k,\ldots,m} \big| \bar X_{t_\ell}^{h,t_k,x,\alpha}\big| \Big\|_p & \leq &
K_p(1+|x|),
\enq
for some positive constant $K_p$, not depending on $(h,t_k,x,i)$.  Moreover, one can also derive the standard estimate for the Euler scheme, as e.g. 
in section 10.2 of  \cite{klopla99}:
\beq \label{estimEuler2}
\sup_{\alpha\in\Ac^h_{t_k,i}} \Big\| \max_{\ell = k,\ldots,m} \big| X_{t_\ell}^{t_k,x,\alpha} -  \bar X_{t_\ell}^{h,t_k,x,\alpha}\big| \Big\|_p 
& \leq & K_p(1+|x|) \sqrt{h}.
\enq
We then associate to the Euler controlled Markov chain, the value functions $\bar{v}_i^h$, $i$ $\in$ $\I_q$, for the optimal switching problem:
\beq
\bar{v}_i^h(t_k,x) &=& \sup_{\alpha\in\Ac_{t_k,i}^h}   \E \Big[ \sum_{\ell=k}^{m-1} f(\bar{X}_{t_\ell}^{h,t_k,x,\alpha},I_{t_\ell}) h  + g(\bar{X}_{t_m}^{h,t_k,x,\alpha},I_{t_m})  \nonumber \\
& & \hspace{3cm}  - \;    \sum_{n=1}^{N(\alpha)} c(\bar{X}_{\tau_n}^{h,t_k,x,\alpha},\iota_{n-1},\iota_n) \Big]. \label{defbarvih}
\enq

The  next result provides the error analysis between $v_i^h$ by $\bar v_i^h$, and thus of the continuous time optimal switching problem $v_i$ by its Euler discrete-time approximation $\bar v_i^h$.

\begin{Theorem} \label{theo2Euler}
There exists a constant $K$ (not depending on $h$) such that  
\beq
\big| v_i^h(t_k,x) - \bar v_i^h(t_k,x) \big| & \leq & K(1+|x|^2) \sqrt{h},
\enq
for all $(t_k,x,i)$ $\in$ $\T_h\times\R^d\times\I_q$.
\end{Theorem}
 
\begin{Remark}
{\rm  The above theorem combined  with Theorem \ref{diff-v-v^h} gives the rate of convergence for the approximation of the 
continuous time optimal switching problem by its Euler discrete-time version:  there exists a positive constant 
$K$ s.t. 
\beq \label{estimdiscret}
|v_i(t_k,x) - \bar v_i^h(t_k,x) | & \leq & K(1 + |x|^{5/2}) \big(h \log(2T/h)\big)^{{1\over 2}}, \;\;\;
\enq 
for  all $(t_k,x,i)$  $\in$ $\T_h\times \R^d\times\I_q$. Moreover if the cost functions $c_{ij}$, $i,i$ $\in$ $\I_q$, do not depend on $x$, then
\beqs 
|v_i(t_k,x) - \bar v_i^h(t_k,x) | & \leq & K(1 + |x|^{2}) h^{{1\over 2}}, \;\;\;
\enqs 
}
\end{Remark}

\vspace{1mm}

\noindent {\bf Proof of Theorem \ref{theo2Euler}.}

\noindent $\bullet$ {\it Step 1.}
For $(t_k,x,i)$ $\in$ $\T_h\times\R^d\times\I_q$, and $\alpha\in\Ac_{t_k,i}^h $ we denote  by
\beqs
J^h(t_k,x,i;\alpha) & = & \E \Big[ \sum_{\ell=k}^{m-1} f(X_{t_\ell}^{t_k,x,\alpha},I_{t_\ell}) h  + g(X_{t_m}^{t_k,x,\alpha},I_{t_m}) 
- \;  \sum_{n=1}^{N(\alpha)} c(X_{\tau_n}^{t_k,x,\alpha},\iota_{n-1},\iota_n) \Big], 
\enqs
so that $v^h_i(t,_k,x) = \sup_{\alpha\in\Ac^h_{t_k,i}}J^h(t_k,x,i,\alpha)$. 
Given  $\alpha\in\Ac_{t_k,i}^h $, let us  define   $F^\alpha_\ell$ $=$ $f(X_{t_\ell}^{t_k,x,\alpha},I_{t_\ell}^{\alpha})$, $c_\ell^\alpha$ $=$ 
$c(X_{t_\ell}^{t_k,x,\alpha},I_{t_{\ell-1}}^{\alpha},I_{t_\ell}^{\alpha})$ and    
$Y^\alpha_\ell$ $=$ $\E\big[\sum_{j=\ell}^m \big(hF^\alpha_j-c_j^\alpha\big)|\Fc_{t_\ell}\big]$,  
for $\ell$ $=$ $k,\ldots,m$.
Consider the stopping time 
\beqs
\tau & = &  \inf \{ t_\ell \geq t_k~:~   J^h(t_\ell,X_{t_\ell}^{t_k,x,\alpha},I_{t_\ell}^\alpha;\alpha^0) \geq Y^{\alpha}_\ell \}\;,
\enqs
where $\alpha^0$ is the strategy with no switches, and 
define $\tilde{\alpha} = (\tilde{\tau}_n, \iota_n)$, with 
\beqs
\tilde{\tau}_n &=& \tau_n \mathbf{1}_{\{\tau_n \leq \tau\}} + \infty \mathbf{1}_{\{\tau_n > \tau\}}.
\enqs
As in the proof of Proposition 2.1, we easily check that 
\beq\label{val tilde alpha> val alpha en k}
Y^{\tilde \alpha}_k &\geq& Y^{ \alpha}_k, 
\enq
 and 
 \beq \label{Y tilde < J0}
 Y^{\tilde \alpha}_\ell  & \geq &  J(t_\ell,X_{t_\ell}^{t_k,x,\tilde{ \alpha}},I_{t_\ell}^{\tilde \alpha};\alpha^0)\;, 
\enq
for all $\ell=k,\ldots,m$.
From  \reff{Y tilde < J0} and the estimates on $X_{t_\ell}^{t_k,x,\alpha}$ in Lemma 
\ref{unifboundXalpha}, we  know that
\beq \label{estimFY}
\E\Big[ \sup_{k\leq \ell \leq m} \big(|Y^{\tilde \alpha}_\ell|^2 + |F^{\tilde \alpha}_\ell|^2 + |c^{\tilde \alpha}_\ell|^2\big) \Big]&\leq& K(1+|x|^2)\;,  
\enq
for some positive constant $K$.  Moreover, by definition, we have:
\beqs \label{discBack}
Y^{\tilde \alpha}_\ell &=& \E\left[Y^{\tilde \alpha}_{\ell+1} | \Fc_{t_\ell} \right] + h F_\ell - c_\ell, \;\;\; \ell = k,\ldots,m-1. 
\enqs
Letting $\Delta M^{\tilde \alpha}_{\ell+1}$ $:=$  $Y^{\tilde \alpha}_{\ell+1} - \E[Y^{\tilde \alpha}_{\ell+1} | \Fc_{t_\ell}]$, we obtain in particular
\beqs
\sum_{\ell=k}^{m-1} c^{\tilde \alpha}_\ell &=& h \sum_{\ell=k}^{m-1} F^{\tilde \alpha}_\ell -  \sum_{\ell=k}^{m-1} \Delta M^{\tilde \alpha}_{\ell+1} + (Y_m^{\tilde \alpha} - Y_k^{\tilde \alpha})\;,
\enqs
and so by \reff{estimFY}
\beq
\E \Big|\sum_{\ell=k}^m c^{\tilde \alpha}_\ell\Big|^2 &\leq& K (1+|x|^2) + 3  \;    \E\left[\left(\sum_{\ell=k}^{m-1} \Delta M^{\tilde \alpha}_{\ell+1}\right)^2 \right] \nonumber \\
&=& K (1+|x|^2) + 3 \;  \E\left[\sum_{\ell=k}^{m-1} |\Delta M_{\ell+1}^{\tilde \alpha}|^2 \right]. \label{ineqc}
\enq
Now by writing that 
\beqs
|Y^{\tilde \alpha}_m|^2 - |Y^{\tilde \alpha}_k|^2 & = &  \sum_{\ell=k}^{m-1} \big(|Y_{\ell+1}^{\tilde \alpha}|^2 - |Y_\ell^{\tilde \alpha}|^2\big) \; = \; \sum_{\ell=k}^{m-1} (Y_{\ell+1}^{\tilde \alpha} - Y_\ell^{\tilde \alpha})(Y_{\ell+1}^{\tilde \alpha} + Y_\ell^{\tilde \alpha} ) \\
&=& \sum_{\ell=k}^{m-1} (\Delta M^{\tilde \alpha}_{\ell+1} - h F^{\tilde \alpha}_\ell + c^{\tilde \alpha}_\ell )(2 Y^{\tilde \alpha}_\ell + \Delta M^{\tilde \alpha}_{\ell+1} - h F^{\tilde \alpha}_\ell + c^{\tilde \alpha}_\ell ),
\enqs
we get
\beqs
\sum_{\ell=k}^{m-1} |\Delta M^{\tilde \alpha}_{\ell+1}|^2 &=& |Y_m^{\tilde \alpha}|^2 - |Y_0^{\tilde \alpha}|^2 - \sum_{\ell=0}^{m-1} h F_\ell^{\tilde \alpha}(h F_\ell^{\tilde \alpha} - 2 Y_\ell^{\tilde \alpha} - 2 c^{\tilde \alpha}_\ell) - 2  \sum_{\ell=0}^{m-1} c_\ell^{\tilde \alpha} Y_\ell^{\tilde \alpha}  \\
& & \;\;\;\;\;\;\; - \;   \sum_{\ell=0}^{m-1} \Delta M_{\ell+1}^{\tilde \alpha}(2 Y_\ell^{\tilde \alpha}- 2 h F_\ell^{\tilde \alpha} + 2 c_\ell^{\tilde \alpha}) - \sum_{\ell=0}^{m-1} |c_\ell^{\tilde \alpha}|^2 .
\enqs
Since $\E\Big[\Delta M^{\tilde \alpha}_{\ell+1} | \Fc_{t_\ell}\Big] =0$, this shows that
\beq
\E\Big[\sum_{\ell=k}^{m-1} |\Delta M_{\ell+1}^{\tilde \alpha}|^2\Big]  &\leq& \E\Big[|Y_m^{\tilde \alpha}|^2 - \sum_{\ell=0}^{m-1} h F^{\tilde \alpha}_\ell(h F^{\tilde \alpha}_\ell - 2 Y^{\tilde \alpha}_\ell - 2 c^{\tilde \alpha}_\ell )- 2 \sum_{\ell=0}^{m-1} c^{\tilde \alpha}_\ell Y^{\tilde \alpha}_\ell \Big] \nonumber \\
&\leq& K (1+|x|^2) + 2 \E \Big[ \Big|\sum_{\ell=0}^{m-1} c^{\tilde \alpha}_\ell Y^{\tilde \alpha}_\ell \Big| \Big], \label{ineqDelta}
\enq
where we used again \reff{estimFY}. Now since $c_\ell \geq 0$,
\beqs
\E \Big[ \Big|\sum_{\ell=0}^{m-1} c^{\tilde \alpha}_\ell Y^{\tilde \alpha}_\ell \Big| \Big] &\leq& \E \Big[ \Big(\sum_{\ell=0}^{m-1} c_\ell \Big) \sup_{k\leq \ell \leq m-1} |Y_\ell^{\tilde \alpha} | \Big] \\
&\leq& \eps \E\Big[\sum_{\ell=k}^{m-1} |\Delta M^{\tilde \alpha}_{\ell+1}|^2 \Big] + K\big(1+\frac{1}{\eps}\big) (1+|x|^2),
\enqs
for all $\eps>0$, by \reff{estimFY}, \reff{ineqc} and Cauchy-Schwarz inequality. Hence taking $\eps$ small enough and plugging this estimate into \reff{ineqDelta}, we obtain
\beqs
\E\Big[\sum_{\ell=k}^{m-1} |\Delta M^{\tilde \alpha}_{\ell+1}|^2\Big]  &\leq& K(1+|x|^2).
\enqs 
 Using \reff{ineqc} one more time and recalling that $N(\tilde \alpha)$ $\leq$  $\eta\sum_\ell c^{\tilde \alpha}_\ell$ for some $\eta$ $>$ $0$ under the uniformly lower bound condition in {\bf (Hc)}, we  thus obtain
\beq \label{estimN(tildealphah)}
\E\big|N(\tilde \alpha)\big|^2 &\leq& K(1+|x|^2).
\enq
Combining this last inequality with \reff{val tilde alpha> val alpha en k}, we get that the supremum in the definition \reff{defvih} of $v_i^h(t_k,x)$ can be taken over 
$\Ac_{t_k,i}^{h,K}(x)$ $=$ $\big\{ \alpha \in \Ac_{t_k,i}^h \mbox{ s.t. } \E|N(\alpha)|^2 \leq K(1+|x|^2) \big\}$.
Using the same argument  with $\bar X^{t_k,x,\alpha}$ instead of $X^{t_k,x,\alpha}$ and   estimate  \reff{estimEuler} on 
$\big\|\bar{X}_{t_\ell}^{h,t_k,x,\alpha}\big\|_2$
we also get that the supremum in the definition  \reff{defbarvih}  $\bar{v}_i^h(t_k,x)$ can be taken over 
$\Ac_{t_k,i}^{h,K}(x)$.
%
%
 
\vspace{1mm}

\noindent $\bullet$ {\it Step 2.}   
Now, for any  $\alpha \in \Ac_{t_k,i}^{h,K}(x)$, we have under {\bf (Hl)} and by Cauchy-Schwarz inequality
\beq
&& \E \Big[ \sum_{\ell=k}^{m-1} h  \big| f(X_{t_\ell}^{t_k,x,\alpha},I_{t_\ell}) - f(\bar{X}_{t_\ell}^{h,t_k,x,\alpha},I_{t_\ell})\big|  
+  \big| g(X_{t_m}^{t_k,x,\alpha},I_{t_m}) - g(\bar{X}_{t_m}^{h,t_k,x,\alpha},I_{t_m}) \big|  \nonumber \\ 
& & \;\;\;\;  + \;   \sum_{n= 1}^{N(\alpha)}  \big|c(X_{\tau_n}^{t_k,x,\alpha},\iota_{n-1},\iota_n) - 
c(\bar{X}_{\tau_n}^{h,t_k,x,\alpha},\iota_{n-1},\iota_n) \big| \Big] \nonumber \\
&\leq& K  \E\Big[ (1 + N(\alpha)) \big( \sup_{k\leq \ell \leq m} \big|X_{t_\ell}^{t_k,x,\alpha} - \bar{X}_{t_\ell}^{h,t_k,x,\alpha}\big| \big)\Big] \nonumber \\
& \leq & K (1 + |x|) \Big\|  \sup_{k \leq \ell \leq m}  \big|X_{t_\ell}^{t_k,x,\alpha} - \bar{X}_{t_\ell}^{h,t_k,x,\alpha}\big| \Big\|_2   \nonumber \\
& \leq &  K (1 + |x|^2) \sqrt{h}, \label{diffvi}
\enq
by \reff{estimEuler2}. Taking the supremum over $\alpha \in \Ac_{t_k,i}^{h,K}(x)$ into 
\reff{diffvi}, this shows that
\beqs
\big| v_i^h(t_k,x) - \bar{v}_i^h(t_k,x) \big| &\leq& K (1+|x|^2)  \sqrt{h}. \label{ineqv-vbar}
\enqs 
\ep

 \section{Approximation schemes by optimal quantization}

\setcounter{equation}{0} \setcounter{Assumption}{0}
\setcounter{Theorem}{0} \setcounter{Proposition}{0}
\setcounter{Corollary}{0} \setcounter{Lemma}{0}
\setcounter{Definition}{0} \setcounter{Remark}{0}

In this section,  for a fixed time discretization step $h$,  we focus on a computational appro\-ximation for the value functions $\bar v_i^h$, $i$ $\in$ 
$\I_q$, defined in \reff{defbarvih}.  To alleviate notations, we shall often omit the dependence on $h$ in the superscripts, and write e.g. 
$\bar v_i$ $=$ $\bar v_i^h$. The corresponding dynamic programming relation for $\bar v_i$  is written in the backward induction:
\beqs
\bar v_i(t_m,x) &=& g_i(x), \\
\bar v_i(t_k,x) &=& \max \Big\{ \E \big[ \bar v_i(t_{k+1},\bar X_{t_{k+1}}^{t_k,x,i}) \big]  + f_i(x) h \; , \; \max_{j\neq i} [ \bar v_j(t_k,x) - c_{ij}(x) ] \Big\},
\enqs
for $k$ $=$ $0,\ldots,m-1$, $(i,x)$ $\in$ $\I_q\times\R^d$, where $\bar X^{t_k,x,i}$ is the solution to the Euler scheme:
\beqs
\bar X_{t_{k +1 }}^{t_k,x,i} & = & F_i^h(x,\vartheta_{k+1}) \; := \; x + b_i(x)h + \sigma_i(x) \sqrt{h} \; \vartheta_{k+1}.
\enqs
Observe  that under the triangular condition  on the switching costs $c_{ij}$ in {\bf(Hc)},  these backward relations  can be written as an  explicit  discrete-time scheme. Indeed, if $\bar v_i(t_k,x) = \bar v_j(t_k,x) - c_{ij}(x)$ for some $j \neq i$, we then have for $l \neq i,j$, 
\beqs
\bar v_j(t_k,x) - c_{ij}(x) \; = \; \bar v_i(t_k,x)  &\geq& \bar v_l(t_k,x) - c_{il}(x) \\
&>& \bar v_l(t_k,x) - c_{ij}(x) - c_{jl}(x),
\enqs
so that $\bar v_j(t_k,x) > \bar v_l(t_k,x) - c_{jl}(x)$. By positivity of the switching costs, we also have 
$$\bar v_j(t_k,x) = \bar v_i(t_k,x) + c_{ij}(x) > \bar v_i(t_k,x) - c_{ji}(x).$$
 It follows that
\beqs
\bar v_j(t_k,x) &=& \E \big[ \bar v_j(t_{k+1},\bar X_{t_{k+1}}^{t_k,x,j}) \big]  + f_j(x) h,
\enqs
and (recalling that $c_{ii}(\cdot)=0$), the backward induction may be rewritten as
\beq
\bar v_i(t_m,x) &=& g_i(x) \label{bacbarvitm2} \\
\bar v_i(t_k,x) &=& \max_{j\in\I_q} \Big\{ \E \big[ \bar v_j(t_{k+1},\bar X_{t_{k+1}}^{t_k,x,j}) \big]  + f_j(x) h - c_{ij}(x) \Big\},
\label{bacbarvitk2}
\enq
for $k$ $=$ $0,\ldots,m-1$, $(i,x)$ $\in$ $\I_q\times\R^d$.   Next, the practical implementation for this scheme requires a computational 
approximation of the expectations arising in the above dynamic programming  formulae, and a  space discretization for the state process 
$X$ valued in $\R^d$.   We shall propose two numerical approximations schemes by optimal quantization methods, the second one in the 
particular case where the state process $X$ is  not controlled by the switching control.

\subsection{A Markovian quantization method}

Let $\X$ be a bounded lattice grid on $\R^d$ with step $\delta/d$ and size $R$, namely $\X$ $=$ $(\delta/d)\Z^d$ $\cap$ 
$B(0,R)$ $=$ $\{ x \in \R^d: x = (\delta/d) z \mbox{ for some } z \in \Z^d, \mbox{ and } |x|\leq \R\}$. 
We then denote by ${\rm Proj}_{\X}$ the projection on the grid $\X$ according to the closest neighbour  rule, which satisfies
\beq \label{proj}
| x- {\rm Proj}_{\X}(x)| & \leq & \max( |x| - R,0) + \delta, \;\;\; \forall x \in \R^d. 
\enq
At each time step $t_k$ $\in$ $\T_h$, and point space-grid $x$ $\in$ $\X$, we  have to compute in \reff{bacbarvitk2} 
expectations in the form $\E \Big[ \varphi(\bar X_{t_{k+1}}^{t_k,x,i})\Big]$, for  $\varphi(.)$ $=$ $\bar v_i^h(t_{k+1},.)$, $i$ $\in$ $\I_q$.   
We shall then use an optimal quantization for the Gaussian random variable $\vartheta_{k+1}$, which consists in approxima\-ting the 
distribution of $\vartheta$ $\leadsto$ $\Nc(0,I_d)$ by the discrete law of a random variable $\hat\vartheta$ of support $N$ points $w_l$, $l$ $=$ $1,\ldots,N$, in $\R^d$, and defined as the projection of $\vartheta$ on the grid $\{w_1,\ldots,w_N\}$ following the closest neighbor rule.  
The grid $\{w_1,\ldots,w_N\}$ is optimized in order to minimize the distorsion error, i.e. the quadratic $L^2$-norm 
$\big\|\vartheta-\hat\vartheta\big\|_2$. This optimal grid and the associated weights 
$\{\pi_1,\ldots,\pi_N\}$ are downloaded from the website: ``http://www.quantize.maths-fi.com/downloads".  We refer to the survey article \cite{pagphapri04}  for more details on the theoretical and computational aspects of optimal quantization methods. In the vein of \cite{pagphapri04b}, we  introduce the quantized Euler scheme:
\beqs
\hat X_{t_{k+1}}^{t_k,x,i} &=& {\rm Proj}_{\X}(F_i^h(x,\hat\vartheta)), 
\enqs
and define the value functions $\hat v_i$ on $\T_m\times\X$, $i$ $\in$ $\I_q$ in backward induction by 
\beqs
\hat v_i(t_m,x) &=& g_i(x)  \\
\hat v_i(t_k,x) &=& \max_{j\in\I_q} \Big\{ \E \big[ \hat v_j(t_{k+1},\hat X_{t_{k+1}}^{t_k,x,j}) \big]  + f_j(x) h - c_{ij}(x) \Big\}, \; k=0,\ldots,m-1. 
\enqs
This numerical scheme   can be computed explicitly according to the following recursive algorithm: 
 \beqs
 \hat v_i(t_m,x) &=& g_i(x), \;\;\; (x,i) \in \X\times\I_q \\
 \hat v_i(t_k,x) &=&  \max_{j\in\I_q}  \Big[ \sum_{l=1}^N \pi_l \; \hat v_j\big(t_{k+1}, {\rm Proj}_{\X}(F_j^h(x,w_l)) \big) + f_j(x)h - c_{ij}(x) \Big], \;\;\; (x,i) \in \X\times\I_q,
 \enqs
for $k$ $=$ $0,\ldots,m-1$.  At each time step, we need to make $O(N)$ computations for each point of the grid $\X$. Therefore, the global complexity of the algorithm 
is of order $O(mN(R/\delta)^d)$.

 \vspace{2mm}

The main result of this paragraph  is to provide an error analysis and rate of convergence for the approximation of 
$\bar v_i$ by $\hat v_i$. 

\begin{Theorem} \label{propv-vhat}There exists a constant $K$ (not depending on $h$) such that
\beqs
\big| \bar v_i(t_k,x) - \hat{v}_i(t_k,x) \big| &\leq&  K \exp\big(K h^{-1} \big\| \vartheta-\hat\vartheta\big\|_2^2\big) \Big(1 + |x|+ \frac{\delta}{h}  \Big) \\
& & \;\;\;\;\;  \Big[ \frac{\delta}{h} +  h^{-1/2} \big\| \vartheta-\hat\vartheta\big\|_2 \Big(1+|x|+ \frac{\delta}{h}\Big) \nonumber \\
&& \;\;\;\;\;\;\;\;   + \;  \frac{1}{Rh}\exp\big(Kh^{-2}\big\| \vartheta-\hat\vartheta\big\|_4^4\big) \Big(1+|x|^2+(\frac{\delta}{h})^2\Big) \Big], 
\nonumber 
\enqs
for all $(t_k,x,i)$ $\in$ $\T_h\times\X\times\I_q$. In the case where the switching costs $c_{ij}$ do not depend on $x$, the above estimation is stengthened into: 
\beqs
\big| \bar v_i(t_k,x) - \hat{v}_i(t_k,x) \big| &\leq&  K \Big[ h^{-1/2} \big\|\vartheta-\hat\vartheta\big\|_2 \exp\big(K h^{-1} \big\| \vartheta-\hat\vartheta\big\|_2^2\big) 
\Big(1 + |x|+ \frac{\delta}{h}   \Big)    \nonumber\\
&& \hspace{3mm}  + \;  \frac{\delta}{h}+ \frac{1}{Rh}\exp\big(K h^{-2} \big\| \vartheta-\hat\vartheta\big\|_4^4\big) \Big(1 + |x|^2+ \big(\frac{\delta}{h}\big)^2 \Big) \Big]. 
\enqs
\end{Theorem}

\begin{Remark} \label{rmEr1}
{\rm  The estimation in Theorem \ref{propv-vhat}  consists of  error terms related to
\begin{itemize}
\item  the space discretization parameters $\delta$, $R$, which have  to be chosen s.t. $\delta/h$ and $1/Rh$ go to zero. 
\item  the quantization error  $\big\|\vartheta-\hat\vartheta\big\|_p$ of the normal distribution $\Nc(0,I_d)$, which converges to zero at a rate $N^{1\over d}$, where 
$N$ is the number of grid points  chosen s.t. $h^{-1\over 2} N^{-1\over d}$ goes to zero.  
\end{itemize}
By combining with the discrete-time approximation error \reff{estimdiscret}, and by choosing grid parameters $\delta$, $1/R$ of order $h^{3\over 2}$,  and a number of points $N$ of order $1/h^{d}$, 
we see that  the error estimate between the value function of the continuous-time optimal switching problem and its approximation by Markovian quantization is of order 
$h^{1\over 2}$.  With these values of the parameters, we then see that the complexity of this Markovian quantization algorithm  is of order $O(1/h^{4d+1})$. 
}
\end{Remark}

\vspace{2mm}

Let us now focus on  the proof of Theorem \ref{propv-vhat}. First, notice from the dynamic pro\-gramming principle that  the value functions 
$\hat v_i$, $i$ $\in$ $\I_q$, admit the Markov control problem representation: 
\beq
\hat v_i(t_k,x) &=& \sup_{\alpha\in\Ac_{t_k,i}^h}   \E \Big[ \sum_{\ell=k}^{m-1} f(\hat X_{t_\ell}^{t_k,x,\alpha},I_{t_\ell}) h  
+ g(\hat X_{t_m}^{t_k,x,\alpha},I_{t_m})  \nonumber  \\
& & \hspace{3cm} -   \sum_{n=1}^{N(\alpha)} c(\hat X_{\tau_n}^{t_k,x,\alpha},\iota_{n-1},\iota_n) \Big], \label{defvih-hat}
\enq
where  $\hat X^{t_k,x,\alpha}$ is defined by
\beqs
\hat X_{t_k}^{t_k,x,\alpha} \; =\;  x, & & \hat X_{t_{\ell+1}}^{t_k,x,\alpha} \; = \;  
{\rm Proj}_{\X}\big(F_{I_{t_\ell}}^h(\hat X_{t_\ell}^{t_k,x,\alpha},\hat\vartheta_{\ell+1})\big), \;\;\;\; k \leq \ell \leq m-1,
\enqs
for $\alpha$ $\in$ $\Ac_{t_k,i}^h$,  and  $\hat\vartheta_{k+1}$, $k$ $=$ $0,\ldots,m-1$, are iid, $\hat\vartheta$-distributed, and independent of 
$\Fc_{t_k}$.  We first prove several estimates on $\hat X^{t_k,x,\alpha}$.

\begin{Lemma} \label{leminter}
For each $p\geq 1$ there exists a constant $K_p$ (not depending on $h$) such that 
\beq
\sup_{\alpha \in \Ac_{t_k,i}^h, k \leq \ell \leq m} \Big\|\hat{X}_{t_\ell}^{t_k,x,\alpha} \Big\|_p &+& 
\sup_{\alpha \in \Ac_{t_k,i}^h, k \leq \ell \leq m-1} \Big\|F^h_{I_{t_\ell}}\big(\hat{X}_{t_\ell}^{t_k,x,\alpha}, \hat\vartheta_{k+1} \big) \Big\|_p  \label{ineqXhat}\\
 &\leq& K_p \exp\left(K_p h^{-p/2} \big\| \vartheta-\hat\vartheta\big\|_p^p\right) \left(1 + |x|+ \frac{\delta}{h}   \right), \nonumber 
\enq
for all $(t_k,x,i)$ $\in$ $\T_h\times\X\times\I_q$.
\end{Lemma}

\ni{\bf Proof.} We fix $(t_k,x,i)$ $\in$ $\T_h\times\X\times\I_q$, $\alpha \in \Ac_{t_k,i}^h$, and denote $\hat X_{t_\ell}$ $=$ $\hat X_{t_\ell}^{t_k,x,\alpha}$, $k \leq \ell \leq m$. Denoting by $\E_l$ the conditional expectation w.r.t. $\Fc_{t_\ell}$, by a standard use of Gronwall's lemma and linear growth of $b_i$, $\sigma_i$, we have
\beq
\E_\ell \Big|F^h_{I_{t_\ell}}(\hat{X}_{t_\ell}, \vartheta_{\ell+1})\Big|^p &\leq& e^{K_p h}\Big|\hat{X}_{t_\ell} \Big|^p + K_p h. \label{ineqhatX1}
\enq
We will use the following convexity inequality : for $a$, $b$ $\in$ $\R_+$, $h$ $\in$ $[0,1]$,
\beq \label{ineqconv}
(a+hb)^p &\leq& (1+K_p h) a^p + K_p h b^p.
\enq
By definition of $F^h$, and the fact that $|{\rm Proj}_{\X}(y)|$ $\leq$ $|y| + \delta$ for all $y$ $\in$ $\R^d$,
\beqs
\Big|\hat{X}_{t_{\ell+1}}\Big| &\leq& \Big|F^h_{I_{t_\ell}}(\hat{X}_{t_\ell}, \vartheta_{\ell+1})\Big| + h^{1/2} \sigma_{I_{t_\ell}}(\hat X_{t_\ell}) \big|\hat\vartheta_{\ell+1}-\vartheta_{\ell+1}\big| + \delta \\
&=& \Big|F^h_{I_{t_\ell}}(\hat{X}_{t_\ell}, \vartheta_{\ell+1})\Big| + h \left(\frac{\sigma_{I_{t_\ell}}(\hat X_{t_\ell}) \big|\hat\vartheta_{\ell+1}-\vartheta_{\ell+1}\big|}{h^{1/2}} + \frac{\delta}{h}\right) 
\enqs
Combining this last inequality with \reff{ineqhatX1}, \reff{ineqconv}, linear growth of $\sigma_i$ and the fact that $\hat\vartheta_{\ell+1},\vartheta_{\ell+1}$ are independent of $\Fc_{t_\ell}$, we obtain
\beqs
\E_\ell\Big|\hat{X}_{t_{\ell+1}} \Big|^p &\leq& (1+K_p h)\big( e^{K_ph} \big|\hat X_{t_\ell}\big|^p + K_p h\big) + K_p h \left(\frac{\sigma_{I_{t_\ell}}(\hat X_{t_\ell})\big\|\vartheta-\hat\vartheta\big\|_p^p}{h^{p/2}} + \frac{\delta^p}{h^p}\right) \\
&\leq& \Big(1+K_p h + K_p h^{1-p/2} \big\|\vartheta-\hat\vartheta\big\|_p^p\Big)\big|\hat X_{t_\ell}\big|^p + K_p h \Big(1+ \big\|\vartheta-\hat\vartheta\big\|_p^p h^{-p/2} + \frac{\delta^p}{h^p}\Big).
\enqs
By induction, taking the expectation, recalling that $h=\frac{T}{m}$, and since $\left( 1 + \frac{y}{m}\right)^m \leq e^y$ for all $y\geq 0$, we obtain
\beqs
\E\Big|\hat{X}_{t_{\ell+1}} \Big|^p &\leq& K_p \exp\left(K_p h^{-p/2} \big\| \vartheta-\hat\vartheta\big\|_p^p\right) \left(1 + |x|^p+ \frac{\delta^p}{h^p} +h^{-p/2}\big\|\vartheta-\hat\vartheta\big\|_p^p   \right) \\
&\leq& K_p \exp\left(K'_p h^{-p/2} \big\| \vartheta-\hat\vartheta\big\|_p^p\right) \left(1 + |x|^p+ \frac{\delta^p}{h^p} \right),
\enqs
for all $k \leq \ell \leq m$. 
The estimate for $F^h(\hat{X}_{t_\ell}, \vartheta_{\ell+1})$ then follows from \reff{ineqhatX1}.
\ep

\vspace{1mm}

\begin{Lemma} \label{lemSup}
There exists some constant $K$ (not depending on $h$) such that
\beq
& & \sup_{\alpha \in \Ac_{t_k,i}^h} \Big\| \sup_{k\leq \ell \leq m} \big|\hat{X}_{t_\ell}^{t_k,x,\alpha} - \bar{X}_{t_\ell}^{t_k,x,\alpha} \big| \Big\|_2  \nonumber  \\
&\leq& K \Big[ h^{-1/2} \big\|\vartheta-\hat\vartheta\big\|_2 \exp\big(K h^{-1/2} \big\| \vartheta-\hat\vartheta\big\|_2\big) \Big(1 + |x|+ \frac{\delta}{h}   \Big)    \nonumber\\
&&  \;\;\;\;\; + \;  \frac{\delta}{h}+ \frac{1}{Rh}\exp\big(K h^{-2} \big\| \vartheta-\hat\vartheta\big\|_4^4\big) \Big(1 + |x|^2+ \big(\frac{\delta}{h}\big)^2 \Big) \Big],  \label{ineqSupXhat-bar}
\enq
for all $(t_k,x,i)$ $\in$ $\T_h\times\X\times\I_q$.
\end{Lemma}
\ni{\bf Proof.} As before we fix $(t_k,x,i)$, $\alpha$ and omit the dependence on $(t_k,x,i,\alpha)$ in $\hat X_{t_\ell}$. Let us first show an estimate on $\Big\|\hat{X}_{t_{\ell+1}} - \bar{X}_{t_{\ell+1}} \Big\|_2$. For $k\leq \ell \leq m-1$, we get
\beq
\Big\|\hat{X}_{t_{\ell+1}} - \bar{X}_{t_{\ell+1}} \Big\|_2 &\leq& \Big\|\hat{X}_{t_{\ell+1}} - F^h_{I_{t_\ell}}(\hat{X}_{t_\ell},\hat \vartheta_{\ell+1}) \Big\|_2 + \Big\|F^h_{I_{t_\ell}}(\hat{X}_{t_\ell}, \hat \vartheta_{\ell+1}) - F^h_{I_{t_\ell}}(\hat{X}_{t_\ell}, \vartheta_{\ell+1}) \Big\|_2 \nonumber \\
&& \;\;\;\; + \;  \Big\|F^h_{I_{t_\ell}}(\hat{X}_{t_\ell}, \vartheta_{\ell+1}) - F^h_{I_{t_\ell}}(\bar{X}_{t_\ell}, \vartheta_{\ell+1}) \Big\|_2.  \label{ineqhatX-barX}
\enq
On the other hand, since
\beqs
\big|y - {\rm Proj}_{\X}(y) \big| &\leq& \delta + |y| {\bf 1}_{\{|y| \geq R\}}  \; \leq \;  \delta + \frac{|y|^2}{R},
\enqs
by inequality \reff{proj}, we have
\beq
\Big\|\hat{X}_{t_{\ell+1}} - F^h_{I_{t_\ell}}(\hat{X}_{t_\ell},\hat \vartheta_{\ell+1}) \Big\|_2 &\leq& \delta + \frac{\Big\|F^h_{I_{t_\ell}}(\hat{X}_{t_\ell},\hat \vartheta_{\ell+1}) \Big\|_4^2}{R}. \label{ineqProj}
\enq
Furthermore by standard estimates for the Euler scheme (see e.g. Lemma A.1 in \cite{pagphapri04b}), we have
\beqs
\Big\|F^h_{I_{t_\ell}}(\hat{X}_{t_\ell}, \vartheta_{\ell+1}) - F^h_{I_{t_\ell}}(\bar{X}_{t_\ell}, \vartheta_{\ell+1}) \Big\|_2 &\leq& (1+K h) \Big\|\hat{X}_{t_\ell} - \bar{X}_{t_\ell} \Big\|_2,
\enqs
and by the linear growth property of $\sigma$ and the fact that $\hat\vartheta_{\ell+1},\vartheta_{\ell+1}$ are independent of $\Fc_{t_\ell}$, 
\beq
\Big\|F^h_{I_{t_\ell}}(\hat{X}_{t_\ell}, \vartheta_{\ell+1}) -F^h_{I_{t_\ell}}(\hat{X}_{t_\ell}, \hat\vartheta_{\ell+1})\Big\|_2 &\leq& K h^{1/2}\left(1 +  \Big\|\hat{X}_{t_\ell} \Big\|_2\right) \big\|\vartheta-\hat\vartheta\big\|_2. \label{ineqhatX2}
\enq
Plugging these three inequalities into \reff{ineqhatX-barX}, we get :
\beqs
\Big\|\hat{X}_{t_{\ell+1}} - \bar{X}_{t_{\ell+1}} \Big\|_2 &\leq& (1+K h) \Big\|\hat{X}_{t_{\ell}} - \bar{X}_{t_{\ell}}\Big\|_2 + K h^{1/2} \left(\Big\|\hat{X}_{t_\ell} \Big\|_2 +1\right) \big\|\vartheta-\hat\vartheta\big\|_2 \\
&& \;\;\;\;\;\;\;\;\;+ \;\delta \;+ \; \frac{\Big\|F^h_{I_{t_\ell}}(\hat{X}_{t_\ell},\hat \vartheta_{\ell+1}) \Big\|_4^2}{R}.
\enqs
Finally since $\hat{X}_{t_{k}}=\bar{X}_{t_k}=x$, we obtain by induction, and using the estimates \reff{ineqXhat} on $\Big\|F^h_{I_{t_\ell}}(\hat{X}_{t_\ell},\hat \vartheta_{\ell+1})\Big\|_4$: 
\beq
\Big\|\hat{X}_{t_{\ell}} - \bar{X}_{t_{\ell}}\Big\|_2 &\leq& K \Big[ h^{-1/2} \big\|\vartheta-\hat\vartheta\big\|_2 \exp\big(K h^{-1} \big\| \vartheta-\hat\vartheta\big\|_2^2\big) 
\Big(1 + |x|+ \frac{\delta}{h}   \Big)   + \frac{\delta}{h}  \nonumber\\
&&  \;\;\;\;\;  + \;  \frac{1}{Rh}\exp\big(K h^{-2} \big\| \vartheta-\hat\vartheta\big\|_4^4\big) \Big(1 + |x|^2+ \big(\frac{\delta}{h}\big)^2  \Big) \Big], \label{ineqXhat-bar}
\enq
for all $k \leq \ell \leq m$. Now by definition of  $\hat{X}_{t_k}$,  $\bar{X}_{t_k}$, we may write for  $k\leq \ell \leq m-1$:
\beqs
\hat{X}_{t_{\ell+1}} - \bar{X}_{t_{\ell+1}} &=& (\hat{X}_{t_\ell} - \bar{X}_{t_\ell}) + h\big( b(\hat{X}_{t_\ell}, I_{t_\ell})-b(\bar{X}_{t_\ell}, I_{t_\ell}) \big) \\
& &  + \;  \sqrt{h} \big( \sigma(\hat{X}_{t_\ell}, I_{t_\ell}) \hat\vartheta_{\ell+1} - \sigma(\bar{X}_{t_\ell}, I_{t_\ell}) \vartheta_{\ell+1} \big) \\
&& + \;  {\rm Proj}_{\X}\big(F^h_{I_{t_\ell}}\big(\hat{X}_{t_\ell}, \hat\vartheta_{\ell+1} )\big) 
- F^h_{I_{t_\ell}} \big(\hat{X}_{t_\ell}, \hat\vartheta_{\ell+1} \big),
\enqs
Since $\hat X_{t_k}$ $=$ $\bar X_{t_k}$ ($=$ $x$), we obtain by induction:
\beq
\left\| \sup_{k\leq \ell \leq m} \left|\hat{X}_{t_\ell} - \bar{X}_{t_\ell} \right| \right\|_2 &\leq& h \sum_{\ell=k}^{m-1} \left\|b(\hat{X}_{t_\ell}, I_{t_\ell})-b(\bar{X}_{t_\ell}, I_{t_\ell}) \right\|_2 \nonumber \\ 
&& + \sqrt{h} \Big\| \sup_{k\leq \ell \leq m} \big| \sum_{r\leq \ell} \sigma(\hat{X}_{t_r}, I_{t_r}) \hat\vartheta_{r+1} 
- \sigma(\bar{X}_{t_r}, I_{t_r}) \vartheta_{r+1}  \big| \Big\|_2 \nonumber \\
&& + \sum_{\ell=k}^{m-1} \Big\|{\rm Proj}_{\X}\big(F^h_{I_{t_\ell}}(\hat{X}_{t_\ell}, \hat\vartheta_{\ell+1} ) \big) 
- F^h_{I_{t_\ell}}\big(\hat{X}_{t_\ell}, \hat\vartheta_{\ell+1} \big) \Big\|_2. \label{bound3}
\enq
We now bound each of the three  terms  in the right hand side  of \reff{bound3}.  First, by the Lipschitz property of $b$ and \reff{ineqXhat-bar}, we have
\beqs
& & h \sum_{\ell=k}^{m-1} \big\|b(\hat{X}_{t_\ell}, I_{t_\ell})-b(\bar{X}_{t_\ell}, I_{t_\ell}) \big\|_2 \\ 
&\leq& K \Big[  h^{-1/2} \big\|\vartheta-\hat\vartheta\big\|_2 \exp\big(K h^{-1} \big\| \vartheta-\hat\vartheta\big\|_2^2\big) \Big(1 + |x|+ \frac{\delta}{h}   \Big)      \nonumber\\
&&   \;\;\;\;\; + \;  \frac{\delta}{h}+ \frac{1}{Rh}\exp\big(K h^{-2} \big\| \vartheta-\hat\vartheta\big\|_4^4\big) \Big(1 + |x|^2+ \big(\frac{\delta}{h}\big)^2 \Big) \Big]. 
\enqs
Next, recalling that $\hat\vartheta_{\ell+1}$ is independent of $\Fc_{t_\ell}$, with distribution law $\hat\vartheta$, and since $\hat\vartheta$ 
is an optimal $L^2$-quantizer  of $\vartheta$, it follows that 
$\E[\hat\vartheta_{\ell+1}|\Fc_{t_\ell}]$ $=$ $\E[\hat\vartheta]$ $=$ $\E[\vartheta]$ $=$ $0$. Thus,  the process $(\sum_{r\leq \ell} \sigma(\hat{X}_{t_r}, I_{t_r}) \hat\vartheta_{r+1} - \sigma(\bar{X}_{t_r}, I_{t_r}) \vartheta_{r+1})_\ell$ is a 
$\Fc_{t_\ell}$-martingale, and from Doob's inequality, we have: 
\beqs
& & \Big\| \sup_{k\leq \ell \leq m} \big| \sum_{r\leq \ell} \sigma(\hat{X}_{t_r}, I_{t_r}) \hat\vartheta_{r+1} 
- \sigma(\bar{X}_{t_r}, I_{t_r}) \vartheta_{r+1}  \big| \Big\|_2  \\
&\leq& K \Big(\E\Big[\sum_{\ell=k}^{m-1} \big|\sigma(\hat{X}_{t_\ell}, I_{t_\ell}) \hat\vartheta_{\ell+1} 
- \sigma(\bar{X}_{t_\ell}, I_{t_\ell}) \vartheta_{\ell+1}\big|^2 \Big] \Big)^{1\over 2}. 
\enqs
By writing from the Lipschitz condition on $\sigma_i$ that
\beqs
\big| \sigma(\hat{X}_{t_\ell}, I_{t_\ell}) \hat\vartheta_{\ell+1} - \sigma(\bar{X}_{t_\ell}, I_{t_\ell}) \vartheta_{\ell+1} \big|^2 &\leq& 
K \Big( \big|\hat{X}_{t_\ell} - \bar{X}_{t_\ell} \big|^2 \big|\vartheta_{\ell+1}\big|^2  \\
& & \;\;\;\;\;\;\;  + \;  \big(1+\big|\hat{X}_{t_\ell}\big|^2\big)\big|\vartheta_{\ell+1} - \hat\vartheta_{\ell+1}\big|^2 \Big), 
\enqs
and since $\vartheta_{\ell+1}, \hat\vartheta_{\ell+1}$ are independent of $\Fc_{t_\ell}$, we then obtain
\beqs
&& \sqrt{h} \Big\| \sup_{k\leq \ell \leq m} \big| \sum_{r\leq \ell} \sigma(\hat{X}_{t_r}, I_{t_r}) \hat\vartheta_{r+1} 
- \sigma(\bar{X}_{t_r}, I_{t_r}) \vartheta_{r+1}  \big| \Big\|_2 \\
&\leq& K  \sup_{k\leq \ell \leq m-1} \Big[ \big\|\hat{X}_{t_\ell} - \bar{X}_{t_\ell} \big\|_2 + \big(1+ \big\|\hat{X}_{t_\ell}\big\|_2\big) 
\big\|\vartheta - \hat\vartheta \big\|_2 \Big]  \\
&\leq& K \Big[  h^{-1/2} \big\|\vartheta-\hat\vartheta\big\|_2 \exp\big(K h^{-1} \big\| \vartheta-\hat\vartheta\big\|_2^2\big) \Big(1 + |x|+ \frac{\delta}{h}   \Big)    \nonumber\\
&&  \;\;\;\;\; + \;  \frac{\delta}{h}+ \frac{1}{Rh}\exp\big(K h^{-2} \big\| \vartheta-\hat\vartheta\big\|_4^4\big) \Big(1 + |x|^2+ \big(\frac{\delta}{h}\big)^2  \Big) \Big],
\enqs
where we used the estimates \reff{ineqXhat} and \reff{ineqXhat-bar}.  Finally the third term in \reff{bound3} is bounded as before by \reff{ineqProj}.
\ep

\vspace{3mm}

\ni {\bf Proof of Theorem  \ref{propv-vhat}.} For $(t_k,x,i)$ $\in$ $\T_h\times\X\times\I_q$, we show as in the proof of Theorem \ref{theo2Euler} that we can restrict to strategies $\alpha\in\Ac^h_{t_k,i}$ such that
\beqs
\E\big|N(\alpha)\big|^2 &\leq& K\Big(1+\sup_{ k \leq \ell \leq m} \Big\|\hat{X}_{t_\ell}^{t_k,x,\alpha} \Big\|_2^2\Big)\;,
\enqs
for some constant $K$, not depending on $(t_k,x,i,h)$.  By using the estimation \reff{ineqXhat},  this means  that  the supremum in the representation \reff{defvih} of $\hat v_i(t_k,x)$  can be taken over the subset 
\beqs
\hat\Ac_{t_k,i}^{h,K}(x) &  =   & \Big\{ \alpha \in \Ac_{t_k,i}^h \mbox{ s.t. } 
\E|N(\alpha)|^2 \leq K \exp\left(K h^{-1} \big\| \vartheta-\hat\vartheta\big\|_2^2\right) \left(1 + |x|^2+ \frac{\delta^2}{h^2} \right)\Big\}.
\enqs
Then,  for 
$\alpha$ $\in$ $\hat\Ac_{t_k,i}^{h,K}(x)$, we have under {\bf (Hl)} and by Cauchy-Schwarz inequality
\beq
&& \E \Big[ \sum_{\ell=k}^{m-1} h  \big| f(\bar X_{t_\ell}^{t_k,x,\alpha},I_{t_\ell}) - f(\hat{X}_{t_\ell}^{t_k,x,\alpha},I_{t_\ell})\big|  
+  \big| g(\bar X_{t_m}^{t_k,x,\alpha},I_{t_m}) - g(\hat{X}_{t_m}^{t_k,x,\alpha},I_{t_m}) \big|  \nonumber \\ 
& & \;\;\;\;  + \;   \sum_{n= 1}^{N(\alpha)}  \big|c(\bar X_{\tau_n}^{t_k,x,\alpha},\iota_{n-1},\iota_n) - 
c(\hat{X}_{\tau_n}^{h,t_k,x,\alpha},\iota_{n-1},\iota_n) \big| \Big] \nonumber \\
&\leq& K  \E\Big[ (1 + N(\alpha)) \big( \sup_{k\leq \ell \leq m} \big|\bar X_{t_\ell}^{t_k,x,\alpha} 
- \hat{X}_{t_\ell}^{t_k,x,\alpha}\big| \big)\Big] \nonumber \\
& \leq & K \exp\left(K h^{-1} \big\| \vartheta-\hat\vartheta\big\|_2^2\right) \big(1 + |x|+ \frac{\delta}{h}  \big)
 \Big\|  \sup_{k \leq \ell \leq m}  \big|\bar X_{t_\ell}^{t_k,x,\alpha} - \hat{X}_{t_\ell}^{t_k,x,\alpha}\big| \Big\|_2   \nonumber \\
& \leq & K \exp\big(K h^{-1} \big\| \vartheta-\hat\vartheta\big\|_2^2\big) \Big(1 + |x|+ \frac{\delta}{h}  \Big) 
\Big[ \frac{\delta}{h} +  h^{-1/2} \big\| \vartheta-\hat\vartheta\big\|_2 \Big(1+|x|+ \frac{\delta}{h}\Big)  \nonumber\\
&& \;\;\;\;\; \;\;  + \;  \frac{1}{Rh}\exp\big(Kh^{-2}\big\| \vartheta-\hat\vartheta\big\|_4^4\big)\Big(1+|x|^2+\big(\frac{\delta}{h}\big)^2\Big) \Big] , \label{diffvihat}
\enq
by Lemma \ref{lemSup}. 
Taking the supremum over $\alpha \in \hat\Ac_{t_k,i}^{h,K}(x)$ in the above inequality, we obtain an estimate for 
$|\bar v_i(t_k,x)-\hat v_i(t_k,x)|$ with an upper bound given by the r.h.s. of \reff{diffvihat}, which gives the required result.  

Finally, notice that in the special case where the switching cost functions $c_{ij}$ do not depend on $x$, we have
\beqs
\big|\bar v_i(t_k,x) - \hat{v}_i(t_k,x) \big| &\leq& \sup_{\alpha \in \Ac_{t_k,i}^h} 
\E \Big[ \sum_{\ell=k}^{m-1} h  \big| f(\bar X_{t_\ell}^{t_k,x,\alpha},I_{t_\ell}) - f(\hat {X}_{t_\ell}^{t_k,x,\alpha},I_{t_\ell})\big| \\
& & \;\;\;\;\;\;\;\;\;\;\;\;\;\;\;  + \;  \big| g(\bar X_{t_m}^{t_k,x,\alpha},I_{t_m}) - g(\hat{X}_{t_m}^{t_k,x,\alpha},I_{t_m}) \big| \Big] \\
& \leq&  K  \sup_{\alpha \in \Ac_{t_k,i}^h, k\leq \ell \leq m} \E\big|\bar X_{t_\ell}^{t_k,x,\alpha} - \hat{X}_{t_\ell}^{t_k,x,\alpha}\big|  \\
& \leq &   K \Big[ h^{-1/2} \big\|\vartheta-\hat\vartheta\big\|_2 \exp\big(K h^{-1} \big\| \vartheta-\hat\vartheta\big\|_2^2\big) \Big(1 + |x|+ \frac{\delta}{h}   \Big)    \nonumber\\
&& \hspace{3mm}  + \;  \frac{\delta}{h}+ \frac{1}{Rh}\exp\big(K h^{-2} \big\| \vartheta-\hat\vartheta\big\|_4^4\big) \Big(1 + |x|^2+ \big(\frac{\delta}{h}\big)^2 \Big) \Big],
\enqs
by  the estimate in Lemma \ref{lemSup}. 
\ep

\subsection{Marginal quantization in the uncontrolled diffusion case}

In this paragraph, we consider the special case where the diffusion $X$ is not controlled, i.e. $b_i=b$, $\sigma_i=\sigma$.  
The Euler scheme for $X$, denoted by $\bar X$, is given  by:
\beqs
\bar X_0 \; = \;  X_0, \;\;\; \bar X_{t_{k+1}} & = & F^h(\bar X_{t_k},\vartheta_{k+1})  \\
&  := &  \bar X_{t_k}  + b(\bar X_{t_k}) h +  \sigma(\bar X_{t_k}) \sqrt{h}  \; \vartheta_{k+1}, 
\;\; k=0,\ldots,m-1,  
\enqs
where $\vartheta_{k+1}$ $=$  $(W_{t_{k+1}}-W_{t_k})/\sqrt{h}$, $k$ $=$ $0,\ldots,m-1$,  are iid, $\Nc(0,I_d)$-distributed, independent of $\Fc_{t_k}$.  
Let us recall the well-known estimate: for any $p$ $\geq$ $1$, there exists some $K_p$ s.t.
\beq \label{estimEulercla}
\big\|\bar X_{t_k} \big\|_p &\leq& K_p (1+ \big\|X_0\big\|_p). 
\enq
Notice that  the backward dynamic programming formulae \reff{bacbarvitm2}-\reff{bacbarvitk2} for $\bar v_i$ can be written in this case as:
\beq
\bar v_i(t_m,.) &=& g_i(.), \;\;\; i \in \I_q \nonumber \\
\bar v_i(t_k,.) &=&  \max_{j\in\I_q} [ P^h \bar v_j(t_{k+1},.) + hf_j - c_{ij} ]. \label{recurbarv}
\enq
Here $P^h$ is the probability transition kernel of the Markov chain $\bar X$, given by:
\beq \label{defPh}
P^h \varphi(x) &=& \E\big[\varphi(\bar X_{t_{k+1}}) | \bar X_{t_k}=x\big] \; = \;  \E [\varphi(F^h(x,\vartheta))],
\enq
where $\vartheta$ is $\Nc(0,I_d)$-distributed.  Let us next  consider the family of discrete-time processes  $(\bar Y_{t_k}^{i})_{k=0,\ldots,m}$, 
$i$ $\in$ $\I_q$, defined by: 
\beqs
\bar Y_{t_k}^i &=& \bar v_i(t_k,\bar X_{t_k}), \;\;\; k =0,\ldots,m, \; i \in \I_q. 
\enqs

\begin{Remark}
{\rm By the Markov property of the Euler scheme $\bar X$ w.r.t. $(\Fc_{t_k})_k$, we see that $(\bar Y_{t_k}^{i})_{k=0,\ldots,m}$, 
 $i$ $\in$ $\I_q$, satisfy   the backward induction:
\beqs
\bar Y_{t_m}^{i} &=& g_i(\bar X_{t_m}) \; = \; g_i(\bar X_T), \;\;\; i \in \I_q \\
\bar Y_{t_k}^{i} &=& \max_{j\in\I_q} \Big\{ \E\big[ \bar Y_{t_{k+1}}^{j} \big| \Fc_{t_k} \big] + h  f_j(\bar X_{t_k})  - c_{ij}(\bar X_{t_k}) \Big\}, 
\;\; k =0,\ldots,m-1,  
\enqs
and is represented as 
\beqs \label{disswitchetendu}
\bar Y_{t_k}^{i} &=&   \esssup_{\alpha\in\Ac_{t_k,i}^h}   \E \Big[ \sum_{\ell=k}^{m-1}  f(\bar X_{t_\ell},I_{t_\ell}) h  + g(\bar X_{t_m},I_{t_m}) 
- \sum_{n=1}^{N(\alpha)} c(\bar X_{\tau_n},\iota_{n-1},\iota_n) \Big| \Fc_{t_k} \Big].
\enqs
On the other hand,  the continuous-time optimal switching problem \reff{defvi} admits a 
representation in terms of the following reflected Backward Stochastic Differential Equations (BSDE): 
\beq 
Y_t^i &=& g_i(X_T) + \int_t^T f(X_s) ds - \int_t^T Z_s^idW_s + K_T^i - K_t^i, \;\;\; i \in \I_q,  \;   0 \leq t\leq T,  \nonumber  \\
Y_t^i & \geq & \max_{j\neq i} [Y_t^j - c_{ij}(X_t)] \;\;\; \mbox{ and } \;  \int_0^T \big(Y_t^i - \max_{j\neq i} [Y_t^j - c_{ij}(X_t)] \big)dK_t^i \; = \; 0.  
\label{BSDE} 
\enq
We know from  \cite{djehampop09}, \cite{hutan10} or \cite{hamzha10} that there exists a unique solution $(Y,Z,K)$ $=$ $(Y^i,Z^i,K^i)_{i\in\I_q}$ solution to 
\reff{BSDE} with   $Y$ $\in$ $\Sc^2(\R^q)$, the set of adapted continuous processes valued in $\R^q$ s.t. $\E[\sup_{0\leq t\leq T}|Y_t|^2]$ $<$ $\infty$, $Z$ $\in$ $\Mc^2(\R^q)$, the set of 
predictable processes valued in $\R^q$ s.t. $\E[\int_0^T |Z_t|^2 dt]$ $<$ $\infty$, and $K^i$ $\in$ $\Sc^2(\R)$, $K^i_0$ $=$ $0$, $K^i$ is nondecreasing.  Moreover,  we have
\beqs
Y_t^i &=&  v_i(t,X_t), \;\;\; i \in \I_q,   \\
& = & \esssup_{\alpha\in\Ac_{t,i}}   \E \Big[ \int_t^T f(X_s,I_s) ds + g(X_T,I_T) 
- \sum_{n=1}^{N(\alpha)}  c(X_{\tau_n},\iota_{n-1},\iota_n) \Big| \Fc_t \Big], \;\; 0 \leq t\leq T.
\enqs
}
\end{Remark}

 \vspace{2mm}

We propose now an optimal quantization method in the vein of \cite{balpag03} for optimal stopping problems, for a computational approximation of 
$(\bar Y_{t_k}^{i})_{k=0,\ldots,m}$.  This is based  on results about optimal quantization of  each marginal distribution of 
the Markov chain $(\bar X_{t_k})_{0\leq k\leq m}$.  Let us recall the construction. 
For each time step $k$ $=$ $0,\ldots,m$, we are given a grid $\Gamma_k$ $=$ $\{x_k^1,\ldots,x_k^{N_k}\}$ of $N_k$ points in $\R^d$, and we
define the quantizer $\hat X_k$ $=$ ${\rm Proj}_k(\bar X_{t_k})$ of $\bar X_{t_k}$ 
where ${\rm Proj}_k$ denotes a closest neighbour projection  
on $\Gamma_k$.  For $N_k$ being fixed, the grid $\Gamma_k$ is  said to be $L^p$-optimal if it minimizes the $L^p$-quantization error: 
$\| \bar X_{t_k} - {\rm Proj}_k(\bar X_{t_k})\|_{_p}$.  Optimal grids $\Gamma_k$ are produced by a stochastic recursive algorithm,  called Competitive Learning Vector Quantization (or also Kohonen Algorithm), and relying on Monte-Carlo simulations of $\bar X_{t_k}$, $k$ $=$ $0,\ldots,m$.  We refer to \cite{pagphapri04} for details about the CLVQ algorithm.  We also compute the transition weights
\beqs
\pi_k^{ll'} &=& \P[ \hat X_{k+1} = x_{k+1}^{l'} | \hat X_k = x_k^l]  
\; = \; \frac{\P\big[(\bar X_{t_{k+1}},\bar X_{t_k}) \in C_{l'}(\Gamma_{k+1})\times C_{l}(\Gamma_{k})\big]}
{\P\big[\bar X_{t_k} \in  C_{l}(\Gamma_{k})\big]},
\enqs
where $C_l(\Gamma_k)$ $\subset$ $\{ x \in \R^d : |x - x_k^l | = \min_{y \in \Gamma_k}|x - y|\}$, $l$ $=$ $1,\ldots,N_k$, is a 
Voronoi tesselation of $\Gamma_k$.  These weights can be computed either  during the CLVQ phase, or by a regular Monte-Carlo simulation once the grids $\Gamma_k$ are settled.  The associated discrete  probability transition $\hat P_k$ from $\hat X_k$ to $\hat X_{k+1}$, $k$ $=$ $0,\ldots,m-1$,  is given by: 
\beqs
\hat P_k \varphi(x_k^l) &:=& \sum_{l'=1}^{N_{k+1}} \pi_{k}^{ll'} \varphi(x_{k+1}^{l'}) 
\; = \; \E\big[\varphi(\hat X_{k+1}) \big| \hat X_k = x_k^l\big]. 
\enqs
One then defines by backward induction the sequence of  $\R^q$-valued functions $\hat v_k$ $=$ $(\hat v_k^i)_{i\in\I_q}$ computed explicitly  on 
$\Gamma_k$, $k$ $=$ $0,\ldots,m$,  by the quantization tree algorithm: 
\beq
\hat v_m^i &=& g_i, \;\;\; i \in \I_q,  \nonumber \\
\hat v_k^i &=& \max_{j\in\I_q} \big[ \hat P_k \hat v_{k+1}^j +  h f_j -  c_{ij} \big], \;\;\; \; k=0,\dots,m-1.  \label{recurhatv}
\enq
The  discrete-time processes  $(\bar Y_{t_k}^{i})_{k=0,\ldots,m}$, $i$ $\in$ $\I_q$, are then approximated by the quantized processes 
$(\hat Y_k^i)_{k=0,\ldots,m}$, $i$ $\in$ $\I_q$ defined by 
\beqs
\hat Y_k^i &=& \hat v_k^i(\hat X_k), \;\;\; k=0,\ldots,m, \; i \in \I_q. 
\enqs

\vspace{2mm}

The rest of this section is devoted to the error analysis  between $\bar Y^{i}$ and $\hat Y^i$.  The analysis follows arguments as in \cite{balpag03b} 
for optimal stopping problems, but has to be slightly modified since the functions $\bar v_i(t_k,.)$ are not Lipschitz in general when the switching costs depend on $x$.  Let us introduce  the subset  $LLip(\R^d)$  of measurable functions $\varphi$ on $\R^d$ satisfying: 
\beqs
| \varphi(x) - \varphi(y) | & \leq & K (1 + |x| + |y|)|x-y|, \;\;\;\;\; \forall x,y \in \R^d,
\enqs
for some positive constant $K$, and denote by 
\beqs
[\varphi]_{_{LLip}} &=& \sup_{x,y\in\R^d, x\neq y} \frac{| \varphi(x) - \varphi(y) | }{(1 + |x| + |y|)|x-y|}. 
\enqs

\begin{Lemma} \label{lemLipvbar}
The  functions $\bar v_i(t_k,.)$,  $k$ $=$ $0,\ldots,m$, $i$ $\in$ $\I_q$, lie in $LLip(\R^d)$, and $[\bar v_i(t_k,.)]_{_{LLip}}$ is bounded by a constant not depending on $(k,i,h)$. 
\end{Lemma}
\ni{\bf Proof.}  We set $\bar v_k^i$ $=$ $\bar v_i(t_k,.)$. From the representation \reff{defbarvih}, we have
\beqs
\bar v_k^i(x) &=& \sup_{\alpha\in\Ac_{t_k,i}^h}   \E \Big[ \sum_{\ell=k}^{m-1}  f(\bar X_{t_\ell}^{t_k,x},I_{t_\ell}) h  
+ g(\bar X_{t_m}^{t_k,x},I_{t_m})  - \sum_{n=1}^{N(\alpha)} c(\bar X_{\tau_n}^{t_k,x},\iota_{n-1},\iota_n) \Big],
\enqs
where $\bar X^{t_k,x}$ is the solution to the Euler scheme starting from $x$ at time $t_k$. From  \reff{estimEulercla}, we see, as in the proof of Theorem \ref{theo2Euler},  that in the above representation for $\bar v_k^i(x)$, one can restrict  the supremum to $\Ac_{t_k,i}^{h,K}(x)$ $=$ $\big\{ \alpha \in \Ac_{t_k,i}^h \mbox{ s.t. } \E|N(\alpha)|^2 \leq K(1+|x|^2) \big\}$ for some positive constant $K$ not depending on $(t_k,x,i,h)$. 
Then, as  in the proof of Theorem \ref{propv-vhat}, we have for any $x,y$ $\in$ $\R^d$, and $\alpha$ $\in$ $\Ac_{t_k,i}^{h,K}(x)$ $\cup$ $\Ac_{t_k,i}^{h,K}(y)$,  
\beqs
&& \E \Big[ \sum_{\ell=k}^{m-1} h  \big| f( \bar X_{t_\ell}^{t_k,x},I_{t_\ell}) - f(\bar{X}_{t_\ell}^{t_k,y},I_{t_\ell})\big|  
+ \big| g(\bar X_{t_m}^{t_k,x},I_{t_m}) - g(\bar{X}_{t_m}^{t_k,y},I_{t_m}) \big| \\ 
& & \;\;\;\;\;\;\;  + \;   \sum_{n=1}^{N(\alpha)} \big|c(\bar X_{\tau_n}^{t_k,x},\iota_{n-1},\iota_n) - c(\bar{X}_{\tau_n}^{t_k,x},\iota_{n-1},\iota_n)\big|\Big] \\
&\leq& K  \big(1 + \big\|N(\alpha)\big\|_2 \big)  \Big\| \sup_{k\leq \ell \leq m} \big| \bar X_{t_\ell}^{t_k,x} - \bar{X}_{t_\ell}^{t_k,y} \big| \Big\|_2  \\
& \leq & K(1 + |x|  + |y|) | x - y|,
\enqs
by standard Lipschitz estimates on the Euler scheme. By taking the supremum over $\Ac_{t_k,i}^{h,K}(x)$ $\cup$ $\Ac_{t_k,i}^{h,K}(y)$ in the above inequality, this shows that 
\beqs
|\bar v_k^i(x) - \bar v_k^i(y) | & \leq & K(1 + |x|  + |y|) | x - y|,
\enqs
i.e.  $\bar v_k^i$ $\in$ $LLip(\R^d)$ with $[\bar v_k^i]_{_{LLip}}$ $\leq$ $K$. 
\ep

\vspace{3mm}

The next Lemma shows that the probability transition kernel of the Euler scheme preserves the growth linear Lipschitz property.  
   
\begin{Lemma} \label{lemLipP}
For any $\varphi$ $\in$ $LLip(\R^d)$, the function $P^h\varphi$ also lies in $LLip(\R^d)$, and there exists some constant $K$, not depending on 
$h$, such that 
\beqs
[P^h \varphi]_{_{LLip}}  & \leq & \sqrt{3}(1 + O(h))  [\varphi]_{_{LLip}},  
\enqs
where $O(h)$ denotes any function s.t. $O(h)/h$ is bounded when $h$ goes to zero. 
\end{Lemma}
\ni{\bf Proof.} From \reff{defPh} and Cauchy-Schwarz inequality, we have for any $x,y$ $\in$ $\R^d$:
\beq
& & |P^h \varphi(x)-P^h \varphi(y)| \nonumber \\
&\leq& \Big(\E \big| \varphi(F^h(x,\vartheta)) - \varphi(F^h(y,\vartheta)) \big|^2 \Big)^{1/2} \nonumber \\
&\leq&  [\varphi]_{_{LLip}}  \Big(\E \big| \big(1+|F^h(x,\vartheta)|+|F^h(y,\vartheta)|)^2 \big|F^h(x,\vartheta)-F^h(y,\vartheta)\big|^2 \Big)^{1/2} \nonumber \\
&\leq&   \sqrt{3} [\varphi]_{_{LLip}}  \Big(\E \big[(1+|F^h(x,\vartheta)|^2+|F^h(y,\vartheta)|^2)|F^h(x,\vartheta)-F^h(y,\vartheta)|^2\big] \Big)^{1\over 2}, \label{interPh}
\enq
where we used the relation $(a+b+c)^2$ $\leq$ $3(a^2+b^2+c^2)$. Since $\vartheta$ has a symmetric distribution, we have 
\beqs
&& \E\Big[\big(1+|F^h(x,\vartheta)|^2+|F^h(y,\vartheta)|^2\big)|F^h(x,\vartheta)-F^h(y,\vartheta)|^2\Big] \\
&=& \frac{1}{2}\E\Big[\big(1+|F^h(x,\vartheta)|^2+|F^h(y,\vartheta)|^2\big)|F^h(x,\vartheta)-F^h(y,\vartheta)|^2 \\
&&  \;\;\;\;\;\;  + \;  \big(1+|F^h(x,-\vartheta)|^2+|F^h(y,-\vartheta)|^2\big)|F^h(x,-\vartheta)-F^h(y,-\vartheta)|^2\Big]
\enqs
A  straightforward calculation gives 
\beqs
&& \frac{1}{2} \Big[ \big(1+|F^h(x,\vartheta)|^2+|F^h(y,\vartheta)|^2\big) |F^h(x,\vartheta)-F^h(y,\vartheta)|^2 \\
&&  \;\;\;\;\;\;\; + \;  \big(1+|F^h(x,-\vartheta)|^2+|F^h(y,-\vartheta)|^2\big)|F^h(x,-\vartheta)-F^h(y,-\vartheta)|^2 \Big] \\
&=& \big(1+|x+hb(x)|^2+|y+hb(y)|^2 + h |\sigma(x) \vartheta|^2 + h|\sigma(y)\vartheta|^2 \big) \big|x-y+h(b(x)-b(y))\big|^2 \\
& & + \;  h|(\sigma(x)-\sigma(y))\vartheta|^2\big(|x+hb(x)|^2+|y+hb(y)|^2 \big)  \\
& & + \;   4h \Big[ \big(x+hb(x) | \sigma(x)\vartheta\big) + \big(y+hb(y)|\sigma(y)\vartheta\big) \Big] 
\big(x-y+h(b(x)-b(y)) |(\sigma(x)-\sigma(y))\vartheta\big) \\
& &  + \;   h^2 (|\sigma(x)\vartheta|^2 + |\sigma(y)\vartheta|^2)|(\sigma(x)-\sigma(y))\vartheta|^2.
\enqs
By Lipschitz continuity of $b$ and $\sigma$, and the fact that $\E|\vartheta|^4$ $<$ $\infty$, we deduce that
\beqs
& & \E\Big[(1+|F^h(x,\vartheta)|^2+|F^h(y,\vartheta)|^2)|F^h(x,\vartheta)-F^h(y,\vartheta)|^2\Big] \\
&\leq& (1 + O(h)) (1+|x|^2+|y|^2)|x-y|^2. 
\enqs
Plugging this last inequality into \reff{interPh}   shows the required result. 
\ep

\vspace{3mm}

We now pass to the main result of this section by providing some a priori estimates  for $\| \bar Y_{t_k} - \hat Y_k\|$ in terms 
of the quantization error $\| \bar X_{t_k} - \hat X_k\|$. 

\begin{Theorem} \label{theoquantif2}
There exists some positive constant $K$, not depending on $h$, such that 
\beq \label{eqbarY-hatY}
\max_{i\in\I_q} \big\| \bar{Y}_{t_k}^{i} - \hat{Y}_k^i \big\|_p &\leq& 
K  \sum_{\ell= k}^m (1+\| X_0\|_r  + \|\hat{X}_\ell \|_r) \big\|\bar{X}_{t_{\ell}} - \hat{X}_\ell \big\|_s,
\enq
for any $k$ $=$ $0,\ldots,m$, and $(p,r,s)$ $\in$ $(1,\infty)$ s.t. $\frac{1}{p}= \frac{1}{r} + \frac{1}{s}$.
\end{Theorem}
\ni{\bf Proof.} We set  $\bar v_k^i$ $=$ $\bar v_i(t_k,.)$, and by misuse of notations, we also set  
$\bar Y_k^i$ $=$ $\bar Y_{t_k}^i$ $=$  $\bar v_k^i(\bar X_k)$. 
From the recursive induction \reff{recurbarv} (resp. \reff{recurhatv}) on $\bar v_k^i$ (resp. $\hat v_k^i$),
 and the trivial inequality $|\max_j \bar a_j - \max_j \hat a_j|$ $\leq$ $\max_j|\bar a_j-\hat a_j|$, we have 
for all $i$ $\in$ $\I_q$: 
\beqs
|\bar{Y}_k^i - \hat{Y}_k^i| &=& | \bar v_k^i(\bar X_{t_k}) - \hat v_k^i(\hat X_k)| \\
&\leq& \max_{j\in\I_q} \big| \big[ P^h \bar{v}_{k+1}^j(\bar{X}_{t_k}) + h f_j(\bar{X}_{t_k}) -  c_{ij}(\bar{X}_{t_k}) \big] 
-  \big[ \hat P_k \hat v_{k+1}^j(\hat X_{k}) + h f_j(\hat X_k) -  c_{ij}(\hat X_k) \big] \big| \\
& \leq& \max_{j\in\I_q} \Big[\big| P^h \bar{v}_{k+1}^j(\bar{X}_{t_k}) - \hat P_k \hat v_{k+1}^j(\hat X_{k})  \big| 
+  h \big| f_j(\bar{X}_{t_k}) - f_j(\hat X_k) \big| 
+ \big| c_{ij}(\bar{X}_{t_k})  - c_{ij}(\hat X_k)\big| \Big]\\
& \leq & K\big| \bar{X}_{t_k} - \hat{X}_k\big|  +  \max_{j\in\I_q} \big| P^h \bar{v}_{k+1}^j(\bar{X}_{t_k}) - \hat P_k \hat v_{k+1}^j(\hat X_{k})  \big| 
\enqs
by the Lipschitz property of $f_j$ and $c_{ij}$, and so 
\beq \label{interY}
\max_{i\in\I_q} \Big\|\bar{Y}_k^i - \hat{Y}_k^i\Big\|_p & \leq & 
K\Big\| \bar{X}_{t_k} - \hat{X}_k\Big\|_p  +  \max_{i\in\I_q} \Big\| P^h \bar{v}_{k+1}^i(\bar{X}_{t_k}) - \hat P_k \hat v_{k+1}^i(\hat X_{k}) \Big\|_p 
\enq
Writing $\hat{\E}_k$ for the conditional expectation w.r.t. $\hat X_{k}$,  we have for any $i$ $\in$ $\I_q$
\beqs
& & \big| P^h \bar{v}_{k+1}^i(\bar{X}_{t_k}) - \hat P_k \hat v_{k+1}^i(\hat X_{k})\big| \\
&\leq& \big| P^h \bar{v}_{k+1}^i(\bar{X}_{t_k}) - P^h \bar{v}_{k+1}^i(\hat{X}_{k})\big| 
+ \big| P^h \bar{v}_{k+1}^i(\hat{X}_{k}) - \hat{\E}_k [P^h  v_{k+1}^i(\bar X_{t_k})]\big| \\
& & \;\;\; + \; \big| \hat{\E}_k [P^h \bar{v}_{k+1}^i(\bar{X}_{t_k})] - \hat P_k \hat v_{k+1}^i(\hat X_{k})\big| \\
&=& \big| P^h \bar{v}_{k+1}^i(\bar{X}_{t_k}) - P^h \bar{v}_{k+1}^i(\hat{X}_{k})\big| + \big|\hat{\E}_k [ P^h \bar{v}_{k+1}^i(\hat{X}_{k}) -  
P^h \bar v_{k+1}^i(\bar X_{t_k})]\big| \\
&& \;\;\;  + \;  \big| \hat{\E}_k [\bar{Y}_{k+1}^i - \hat{Y}_{k+1}^i]\big|.
\enqs
Since $\hat{\E}_k$ is a $L^p$-contraction, we then obtain
\beq
& & \Big\| P^h \bar{v}_{k+1}^i(\bar{X}_{t_k}) - \hat P_k \hat v_{k+1}^i(\hat X_{k})\Big\|_p \nonumber \\ 
&\leq& 2 \Big\| P^h \bar{v}_{k+1}^i(\bar{X}_{t_k}) - P^h \bar{v}_{k+1}^i(\hat{X}_{k})\Big\|_p + \Big\|\bar{Y}_{k+1}^i - \hat{Y}_{k+1}^i\Big\|_p\nonumber  \\
&\leq& K(1+O(h))\Big\|\big(1+ \big| \bar{X}_{t_k} \big| + \big| \hat{X}_k \big| \big)\big|\bar{X}_{t_k} - \hat{X}_k\big| \Big\|_p 
+ \Big\|\bar{Y}_{k+1}^i - \hat{Y}_{k+1}^i\Big\|_p \nonumber \\
&\leq& K(1+O(h))\big(1+ \big\|X_0\big\|_r  + \big\|  \hat{X}_k \big\|_r \big) \Big\|\bar{X}_{t_k} - \hat{X}_k\Big\|_s 
+ \Big\|\bar{Y}_{k+1}^i - \hat{Y}_{k+1}^i\Big\|_p, \label{interPh2}
\enq
where we  used Lemmata \ref{lemLipP} and \ref{lemLipvbar}, H\"older's inequality and \reff{estimEulercla}.  
Substituting \reff{interPh2} into \reff{interY}, we get 
\beqs
& & \max_{i\in\I_q} \Big\|\bar{Y}_k^i - \hat{Y}_k^i\Big\|_p  \\
& \leq &   K(1+O(h))\Big(1+ \big\|X_0\big\|_r  + \big\|  \hat{X}_k \big\|_r \Big) \Big\|\bar{X}_{t_k} - \hat{X}_k\Big\|_s 
+  \max_{i\in\I_q}  \Big\|\bar{Y}_{k+1}^i - \hat{Y}_{k+1}^i \Big\|_p,
\enqs
for all $k$ $=$ $0,\ldots,m-1$. Since $\max_{i\in\I_q} \big\|\bar{Y}_m^i - \hat{Y}_m^i\big\|_p$ $=$  
$\max_{i\in\I_q}\big\| g_i(\bar X_{t_m}) - g(\hat X_m)\big\|_p$ $\leq$ $K\big\|\bar X_{t_m} - \hat X_m\big\|_p$ by the Lipschitz condition on $g_i$, we  conclude by induction.  
\ep

\vspace{2mm}

\begin{Remark} \label{rmQerror}
{\rm  Assume that $\hat X_k$ is chosen to be an $L^2$-optimal quantizer of $\bar X_{t_k}$  for each $k$ $=$ $0,\ldots,m$. It is in particular a stationary quantizer in the sense that $\E[\bar X_{t_k}|\hat X_k]$ $=$ $\hat X_k$ (see \cite{pagphapri04}), and by Jensen's inequality, we deduce that  $\big\|\hat X_k\big\|_2$ $\leq$ $\|\bar X_{t_k}\big\|_2$.  Recalling \reff{estimEulercla},  
the inequality  \reff{eqbarY-hatY} in Theorem \ref{theoquantif2}  gives
\beqs
\max_{i\in\I_q} \big\| \bar{Y}_{t_k}^{i} - \hat{Y}_k^i \big\|_1 &\leq& 
K (1+\big\| X_0\big\|_2)   \sum_{\ell= k}^m \big\|\bar{X}_{t_{\ell}} - \hat{X}_\ell \big\|_2,
\enqs
for all $k$ $=$ $0,\ldots,m$.  In particular, if $X_0$ $=$ $x_0$ is deterministic, then $\hat X_0$ $=$ $x_0$, and we have an error estimation  by quantization of the value function function for the discrete-time optimal switching problem at the initial date measured by:
\beq \label{estimvhat}
\max_{i\in\I_q} \big| \bar v_i(0,x_0) - \hat v_0^i(x_0)\big| & \leq & K(1+|x_0|) \sum_{k=1}^m \big\|\bar{X}_{t_{k}} - \hat{X}_k \big\|_2
\enq
Suppose that one has at hand a global stack of $\bar N$ points for the whole space-time grid, to be dispatched with $N_k$ points for each $k$th-time step, i.e. $\sum_{k=1}^m N_k$ $=$ $\bar N$. 
Then, as in \cite{balpag03b},  in the case of uniformly elliptic diffusion with bounded Lipschitz coefficients 
$b$ and $\sigma$, one can optimize over the $N_k$'s by using the rate of convergence for  the miminal 
$L^2$-quantization error given  by Zador's theorem: 
\beqs
\big\|\bar{X}_{t_{k}} - \hat{X}_k \big\|_2 & \sim &  \frac{J_{2,d} \big\| \varphi_k \big\|_{\frac{d}{d+2}}^{1\over 2} } {N_k^{1\over d}} \;\;\;
 \mbox{ as } N_k \rightarrow \infty,
\enqs
where $\varphi_k$ is the probability density function of $\bar X_{t_k}$, and $\big\|\varphi\big\|_r$ $=$ $(\int |\varphi(u)|^r du)^{1\over r}$.  
From \cite{baltal96}, we have the bound $\big\| \varphi_k \big\|_{\frac{d}{d+2}}^{1\over 2}$  $\leq$ $K \sqrt{t_k}$,
for some constant $K$ depending only on $b$, $\sigma$, $T$, $d$.  Substituting into \reff{estimvhat}  with Zador's theorem, we obtain
\beqs
\max_{i\in\I_q} \big| \bar v_i(0,x_0) - \hat v_0^i(x_0)\big| & \leq & K(1+|x_0|) \sum_{k=1}^m  \frac{\sqrt{t_k}}{N_k^{1\over d}}.
\enqs
For  fixed $h$ $=$ $T/m$ and $\bar N$, the sum in the upper bound of the above inequality is minimized  over  the size of the grids $\Gamma_k$, $k$ $=$ 
$1,\ldots,m$ with
\beqs
N_k &=&  \left \lceil \frac{  t_k^{\frac{d}{2(d+1)}}  \bar N }{ \sum_{k=1}^m t_k^{\frac{d}{2(d+1)}}}  \right \rceil,
\enqs
where $\lceil x\rceil :=\min\{k\!\in\N,\; k\ge x\}$, and we have a global rate of convergence given by: 
\beqs
\max_{i\in\I_q} \big| \bar v_i(0,x_0) - \hat v_0^i(x_0)\big| & \leq & \frac{K(1+|x_0|)}{h(\bar Nh)^{1\over d}}. 
\enqs

Actually even with no extra assumptions on $b$ and $\sigma$, we have the same estimate, since for all $r$ $>$ $0$,
\beqs
\big\|\bar{X}_{t_{k}} - \hat{X}_k \big\|_2 \leq C_{2,r} \big\|\bar{X}_{t_{k}}\big\|_{2+r} N_k^{-1/d} \leq K N_k^{-1/d},
\enqs
see Lemma 1 in \cite{luspag08}.

By combining with  the estimate  \reff{estimdiscret},  we obtain an error  bound between the value function of the continuous-time optimal switching problem and its approximation by marginal 
quantization of order $h^{1\over 2}$ when choosing a number of points by grid $\bar Nh$ of order $1/h^{\frac{3d}{2}}$.  This has to be compared with the number of points $N$ of lower  order 
$1/h^{d}$ in the Markovian quantization approach, see Remark \ref{rmEr1}.  The complexity of this marginal quantization algorithm is of order 
$O\left(\sum_{k=1}^m N_k N_{k+1}\right)$. In terms of $h$, if we take $N_k = \bar N h = 1/h^{\frac{3d}{2}}$, we then need $O(1/h^{3d+1})$ operations to compute the value function.  Recall that the Markovian quantization method requires a complexity of higher 
order $O(1/h^{4d+1})$, but  provides  in compensation an approximation of the value function in the whole space grid $\X$.  
}
\end{Remark}

\section{Numerical tests}

\setcounter{equation}{0} \setcounter{Assumption}{0}
\setcounter{Theorem}{0} \setcounter{Proposition}{0}
\setcounter{Corollary}{0} \setcounter{Lemma}{0}
\setcounter{Definition}{0} \setcounter{Remark}{0}

We test our quantization  algorithms  by comparison results with explicit formulae for optimal switching problems derived from chapter 5 in \cite{pha09}. 
The formulae are obtained  for infinite horizon problems, that we adapt to our case by taking as the final gain the (discounted) value function for the infinite horizon problem.

We consider a two-regime switching problem where the diffusion is independent of the regime and follows a geometric Brownian motion, i.e. 
$b(x,i) = b x$, $\sigma(x,i) = \sigma x$, and the switching costs are constant $c(x,i,j) = c_{ij}$ ,$i,j=1,2$. The profit functions are in the form 
$f_i(t,x) = e^{-\beta t} k_i x^{\gamma_i}$, $i$ $=$ $1,2$.  From Theorem 5.3.5 in \cite{pha09}),  the value functions are given by:
\beqs
v_1(0,x) &=& \left\{ \begin{array}{cl} A_1 x^{m^+} + K_1 k_1 x^{\gamma_1}, & \;\;\; x < \underline{x}_1^* \\
B_2 x^{m^-} + K_2 k_2 x^{\gamma_2} - c_{12}, & \;\;\; x \geq \underline{x}_1^* \end{array}\right. \\
v_2(0,x) &=& \left\{ \begin{array}{cl} A_2 x^{m^+} + K_2 k_2 x^{\gamma_2}, & \;\;\; x < \underline{x}_2^* \\
A_1 x^{m^+} + K_1 k_1 x^{\gamma_1} - c_{21} & \;\;\;  \underline{x}_2^* \leq x \leq \overline{x}_2^*  \\
B_2 x^{m^-} + K_2 k_2 x^{\gamma_2}, & \;\;\; x > \overline{x}_2^*\end{array}\right. ,
\enqs
where $A_i$, $B_i$, $K_i$, $\underline{x}_2^*$ and $\overline{x}_2^*$ depend explicitly on the parameters. In the sequel, we take for value of the parameters:  
\beqs
b = 0,  \; \sigma= 1, \;  c_{01} = c_{10} = 0.5, \; k_1=2, k_2 = 1, \; \gamma_1 = 1/3, \; \gamma_2 = 2/3, \; \beta = 1. 
\enqs

We  compute the value function in regime 2 taken at $X_0=3.0$ by means of the first algorithm (Markovian quantization). We take $R= 10 X_0$ and vary $m,\delta$ and $N$. The results are compared with the exact value in Table 1. Notice  that the algorithm seems to be quite robust and provides  good results even when 
$\delta m$ and $\frac{m}{R}$ do not satisfy the constraints given by our theoretical estimates in Remark \ref{rmEr1}. 

In Table 2, we have computed the value with the marginal quantization algorithm. We make vary the number of time steps $m$ and the total number of grid points $\bar N$ (dispatched between the different time steps as described in Remark \ref{rmQerror}). We have used optimal quantization of the Brownian motion, and the transition probabilities $\pi_k^{ll'}$ were computed by Monte-Carlo simulations with $10^6$ sample paths (for an analysis of the error induced by this Monte-Carlo approximation, see Section 4 in \cite{balpag03}). We have also indicated the time spent for these computations. Actually,  almost all of this time comes from the Monte-Carlo computations, as the tree descent algorithm is very fast (less than 1s for all the tested parameters). 

For the two methods, we look at the impact of the quantization number for each time step (resp. $N$ and $\bar N h$) on the precision of the results. As our theoretical estimates showed (see Remarks \ref{rmEr1} and \ref{rmQerror}), for the first method,  increasing $N$ higher than $h^{-1}$ does not seem to improve the precision, whereas for the second method, we can see for several values of $h$ that changing $\bar N h$ from $h^{-1}$ to $h^{-2}$ or $h^{-3}$ improves the precision. 

Comparing the two tables, the first method seems to provide precise estimates with slightly faster computation times, and it has the further advantage of computing simultaneously the value functions at any points of the space discretization grid $\X$. However, since most of the time spent by our second algorithm was devoted to the calculation of  the transition probabilities $\pi_k^{ll'}$, if these were  computed beforehand and stored offline, the marginal quantization method becomes more competitive.

\vspace{3mm}

\begin{table}[h]
\centering
\begin{tabular}{|c|c|c|c|}
\hline
 $(m,1/\delta,N)$ & $\hat{v}_2(0,3.0)$ & Numerical error ($\%$) & Algorithm time (s)\\
\hline
(10,10,10) & 2.1925 & 3.0 & 0.2 \\
(10,10,100) & 2.1863 & 2.7 & 0.5 \\
(10,10,1000) & 2.1852& 2.7 & 1.4 \\
(10,100,1000) & 2.1882 & 2.8 & 8.5 \\
(10,100,5000) & 2.1882 & 2.8 & 40 \\
(100,10,100) & 2.1218 & 0.31 & 1.0 \\
(100,10,1000) & 2.1213 & 0.33 & 8.0 \\
(100,10,5000) & 2.1213 & 0.33 & 39 \\
(100,100,100) & 2.1250 & 0.16 & 8.6\\
(100,100,1000) & 2.1250 & 0.16 & 82 \\
\hline
Exact value & 2.1285 & &  \\
\hline
\end{tabular}
\label{tbl1}
\caption{Results obtained by Markovian quantization}
\end{table}


\begin{table}[h]
\centering
\begin{tabular}{|c|c|c|c|}
\hline
 $(m,\bar N)$ & $\hat{Y}_0^2$ & Numerical error ($\%$) & Algorithm time (s)\\
\hline
(10,100) & 2.2080 & 3.7 & 4.4 \\
(10,1000) & 2.2174 & 4.2 & 4.9\\
(10,10000) & 2.1276 & 0.04 & 5.8\\
(100,1000) & 2.1233 & 0.24 & 36 \\
(100,10000) & 2.1316 & 0.15 & 48 \\
(100,50000) & 2.1301 & 0.07 & 65\\
(1000,10000) & 2.1161 & 0.58 & 353 \\
(1000,50000) & 2.1213 & 0.34 & 498\\
\hline
\end{tabular}
\label{tbl2}
\caption{Results obtained by marginal quantization}
\end{table}

 \end{document}